\documentclass[onecolumn,final,a4paper]{elsart3p}

\def\ElsevierStyle

\usepackage{float}
\usepackage{comment}
\usepackage{amsmath}       %   AMS equation environments
\usepackage{amssymb}       %   AMS symbol fonts (e.g., \boldsymbol{}

\usepackage{chapterbib}    %   Bibtex
\usepackage{color}
\usepackage{comment}
\usepackage{bm,bbm}
\usepackage{overpic}

\usepackage{times}

\usepackage{epsfig}
\usepackage{graphicx,graphics,rotating}
\usepackage{subfigure}
\usepackage{multicol}
\usepackage{multirow}

\usepackage{booktabs,siunitx}

\newcommand{\TABROW}{}
\newcommand{\citeappx}[1]{}

\newcounter{numbs}
\newcounter{numbi}
\newcounter{numbii}

\definecolor{MyDarkGreen}{rgb}{0,0.45,0}%%{0.80,0.20,0.20} 

\def\trait #1 #2 #3 {\vrule width #1pt height #2pt depth #3pt}
\def\fin{\hfill
        \trait .3 5 0
        \trait 5 .3 0
        \kern-5pt
        \trait 5 5 -4.7
        \trait 0.3 5 0
\medskip}

\newcommand{\REAL}{\mathbbm{R}}

\newcommand{\TERM}[2]{\textbf{(#1)}}

\renewcommand{\cv}{\mathbf{c}}

\newcommand{\ds}{d}

\newcommand{\fs}{f}
\newcommand{\gs}{g}

\newcommand{\js}{j}
\newcommand{\ks}{k}
\newcommand{\ls}{l}
\newcommand{\ms}{m}
\newcommand{\ns}{n}

\newcommand{\ps}{p}

\renewcommand{\ss}{s}
\newcommand{\ts}{t}
\newcommand{\us}{u}
\newcommand{\vs}{v}
\newcommand{\ws}{w}
\newcommand{\xs}{x}

\newcommand{\Ds}{D}

\newcommand{\Hs}{H}

\newcommand{\Ks}{K}

\newcommand{\Ms}{M}
\newcommand{\Ns}{N}

\newcommand{\Ts}{T}

\newcommand{\Vs}{V}

\newcommand{\calM}{\mathcal{M}}
\newcommand{\calN}{\mathcal{N}}

\newcommand{\calS}{\mathcal{S}}

\newcommand{\LTWO}  {L^2}

\newcommand{\PS}[1] {\mathbbm{P}_{#1}}

\newcommand{\dx}{\,dx}
\newcommand{\dt}{\,dt}

\newcommand{\pt}{\partial_t}
\newcommand{\px}{\partial_x}

\newcommand{\pxx}{\partial_{xx}}

\newcommand{\VsN}{\Vs_N}
\newcommand{\vsN}{\vs_N}
\newcommand{\usN}{\us_N}
\newcommand{\fsN}{\fs_N}

\newcommand{\lsN}{\ls^N}
\newcommand{\xiN}{\xi^N}
\newcommand{\wsN}{\ws^N}

\newcommand{\fsh}{\widehat{\fs}}

\newcommand{\phiN}{\phi_N}
\newcommand{\wsa}{\ws_{\alpha}}

\newcommand{\vsh}{\widehat{\vs}}
\newcommand{\Vsh}{\widehat{\Vs}}
\newcommand{\ush}{\widehat{\us}}

\newcommand{\Hsp}{\Hs^{\prime}}
\newcommand{\asp}{\alpha^{\prime}}

\newcommand{\Hspp}{\Hs^{\prime\prime}}

\newcommand{\EOD}{\end{document}}

\begin{document}

\begin{frontmatter}

  \title{  A decision-making machine learning approach in Hermite spectral
  approximations of partial differential equations}
  %{Machine learning selection of the scaling coefficient in the
  %  Hermite spectral approximation of partial differential equations }

  \author[Camerino]   {L. Fatone,}
  \author[Unimore,CNR]{D. Funaro}
  \author[CNR]        {and G. Manzini}

  \address[Unimore]{
    Dipartimento di Scienze Chimiche e Geologiche,
    Universit\`a degli Studi di Modena  e Reggio Emilia, Via campi 103, 41125 Modena
    Italy;
    %\emph{e-mail: daniele.funaro@unimore.it}
  }
  \address[CNR]{
    Istituto di Matematica Applicata e Tecnologie Informatiche,
    Consiglio Nazionale delle Ricerche,
    via Ferrata 1,
    27100 Pavia,
    %\emph{e-mail: daniele.funaro@unimore.it}
  }
  \address[Camerino]{
    Dipartimento di Matematica,
    Universit\`a degli Studi di Camerino,
    Via Madonna delle Carceri 9, 62032 Camerino,
    Italy;
    %\emph{e-mail: lorella.fatone@unicam.it}
  }
 
  %\maketitle
  %% \input{abstract.tex}
  %%
  \begin{abstract}
    The accuracy and effectiveness of Hermite spectral methods for the
    numerical discretization of partial differential equations on
    unbounded domains, are strongly affected by the amplitude of the
    Gaussian weight function employed to describe the approximation
    space.
    This is particularly true if the problem is under-resolved, i.e.,
    there are no enough degrees of freedom.
    The issue becomes even more crucial when the equation under study is
    time-dependent, forcing in this way the choice of Hermite functions
    where the corresponding weight depends on time.
    In order to adapt dynamically the approximation space, it is here
    proposed an automatic decision-making process that relies on machine
    learning techniques, such as deep neural networks and support vector
    machines.
    The algorithm is numerically tested with success on a simple 1D
    problem, but the main goal is its exportability in the context of
    more serious applications.
  \end{abstract}
  

\end{frontmatter}

\newcommand{\PIPPO}{\mbox{\textbf{PIPPO}}}
\newcommand{\xsh}{\widehat{x}}
\newcommand{\psh}{\widehat{p}}
\renewcommand{\ush}{\widehat{u}}
\newcommand{\DISP}{}

\raggedbottom

%% PAPER
%% \input{sec1_introduction.tex}
%% \input{sec2_Hermite.tex}
%% \input{sec3_quadrature.tex}
%% \input{sec4_algorithm.tex}
%% \input{sec5_MLgaussiane.tex}
%% \input{sec6_nonhomogeneous_case.tex}
%% \input{sec7_Vlasov_conclusions.tex}

% Hey Emacs, this is -*-latex-*-

%% sec1

\section{Introduction}
\label{sec1:introduction}
The aim of the present paper is to show how Machine Learning (ML) techniques can be employed,
via a suitable implementation, in
the field of spectral methods for PDEs, using Hermite polynomials as approximating functions.
%%
%% Our approach is for the moment academical, but, as it will be
%% explained in the following paragraphs, it is inspired by precise needs
%% emerging in the context of plasma physics.
%%
To ease the exposition, we approach this subject in an ``academic way'' by just considering
the 1D time-dependent heat equation.
However, we are motivated by precise needs emerging in the context of
plasma physics, which we will explain in the  paragraphs to follow.
In the first part of this introduction, we provide a quick historical
background to the role of Hermite functions as 
approximation basis for PDEs, and we postpone to the second part
a brief discussion on the ML algorithms of our interest.

%%\subsection{Background material: Hermite functions in PDEs}
%%
Generalized Hermite functions are constructed as $\Hs_{n}\wsa$, where $\Hs_n$ is
the $n$-th Hermite polynomial and
$\wsa(x)=\exp(-\alpha^2(\xs-\beta)^2)$ is the weight function.
We refer to $\alpha$ and $\beta$ as the scaling and the shifting
factors.
A first convergence analysis of the Hermite spectral approximation of
the 1D heat equation with both Galerkin and collocation methods was
provided in \cite{Funaro-Kavian:1991}, and for the Burgers equation in
\cite{Guo:1999}, by assuming $\alpha=1$ and $\beta=0$.
Hermite approximations of functions and solutions to PDEs may however
perform poorly when there is not a sufficient number of degrees of
freedom, and inappropriate
choices of the above parameters do not really help improving the outcome.
For example, if $\alpha=1$ and $\beta=0$, the approximation of
$\sin\xs$ in the interval $[-N,N]$ requires at least $N^2$ expansion
terms~\cite[pp.~44--45]{Gottlieb-Orszag:1977}.
Despite this very pessimistic anticipation, we known nowadays that
more effective methods can be obtained by choosing $\alpha$ and
$\beta$ carefully.
As pointed out in~\cite{Tang:1993}, we only need a reasonable number
of terms $N$ in the truncated expansion if $\beta=0$ and the scaling
factor $\alpha$ is taken as
$\ max_{1\leq\js\leq\Ns}\big(\xiN_j\big)\slash{N}$, where $\xiN_j$,
$j=1,\ldots,\Ns$ are the roots of $\Hs_{N}$.  Such a significant
improvement has reopened the way towards the design of Hermite
spectral methods that are computationally efficient for solving
physics and engineering problems.
On the other hand,  this  poses the absolutely
nontrivial problem of how to choose  acceptable values for the
scaling and shifting parameters.

More recently, generalized Hermite functions have been applied to
solving time dependent problems, assuming that 
$\alpha=\alpha(t)$ and  $\beta=\beta(t)$ may
also change in time.
Noteworthy developments in this direction are found for instance 
in~\cite{Ma-Sun-Tang:2005,Luo-Yau-Yau:2015,Wang-Zhang-Qiong:2016}.
Those papers
cover a wide range of applications ranging from the discretization of
linear and nonlinear diffusion and convection-diffusion equations, to
the generalized Ginzburg–Landau equation and the Fokker-Planck
equation.

It is important to mention that  Hermite functions provide the most natural framework for
the numerical approximations of the distribution function solving
the Vlasov equation, which is the most basic mathematical model in noncollisional
plasma physics~\cite{Grad:1949}.
Indeed, Hermite functions are directly linked to the Maxwellian
distribution function that describes a noncollisional plasma close to
the equilibrium state.
When the plasma develops a strong non-Maxwellian behavior, an
impractical number of Hermite basis functions could be needed, thus
making a numerical solution too expensive.
In~\cite{Schumer-Holloway:1998}, it was shown through numerical
experiments that, even in bad situations, using weight functions of the
form $\,\exp\big(-\alpha^2(\xs-\beta)^2\big)$ with a careful choice of
 $\alpha$ and  $\beta$, it is possible to reduce
significantly  the computational burden by orders of
magnitude.
This fact was exploited in successive works for investigating the
properties of Vlasov-based plasma
models~\cite{Camporeale-Delzanno-Lapenta-Daughton:2006,Delzanno-Roytershtein:2018} and developing
computationally efficient numerical methods (see, e.g.:
\cite{Fatone-Funaro-Manzini:2020,Fatone-Funaro-Manzini:2019a,Fatone-Funaro-Manzini:2019b,Funaro-Manzini:2020,Manzini-Funaro-Delzanno:2017,Delzanno:2015},
and the code implementations described
in~\cite{Manzini-Delzanno-Vencels-Markidis:2016,Vencels-Delzanno-Manzini-Markidis-BoPeng-Roytershteyn:2015}).

To the best of our knowledge, however, the selection of reasonable
values for $\alpha$ and $\beta$ (also depending on time) is still an
open problem, since only partial, often unsatisfactory, answers have
been given, sometimes limited to specific coefficients that
must be somehow imposed by the user at the beginning of the
simulation.
Here, our contribution is finalized to the automatic detection of the most appropriate value
of $\alpha$, whereas, for simplicity $\beta$ will remain equal to zero. 
The methodology is based on a ML approach that tries to recognize the
shape of the numerical solution while it evolves in time, and select
dynamically the proper value of $\alpha$.
The ML algorithm is initially trained by means of a suitable
input set, where each one of its elements brings along the correct value of
$\alpha$.
As the numerical solution advances in time, the so
instructed ML algorithm evaluates periodically the most reasonable value of
$\alpha$ at that time. The computation then continues
with the updated value.

%% MACHINE LEARNING ALGORITHMS, neural networks
ML algorithms are nowadays widely applied in the
numerical treatment of differential equations arising in science and
engineering.
Many of such learning algorithms have been designed in recent times
for solving classification and regression problems, in the purpose of
achieving better effectiveness and efficiency in tackling physics and
engineering applications with
computers~\cite{Karniadakis-Hesthaven:2019}.
A review of these developments is beyond the scope
of this article, so we briefly recall some background material useful
for the discussion of our problem, without any claim of being
exhaustive.
Neural Networks (NN) can be used to approximate sufficiently regular
functions and their derivatives, as well as solutions to differential equations
~\cite{Cybenko:1989,Pinkus:1999,NguyenThien-TranCong:1999}.
Moreover, NN can reduce the costs of modelling
computational domains with complex geometries when solving forward and
inverse problems~\cite{Ramuhalli-Upda-Upda:2005}, and can be used as a substitute
to conventional constitutive material
models~\cite{Stoffel-Bamer-Markert:2018}.
A detailed introduction to these topics can be found
in~\cite{Yadav-Yadav-Kumar:2015}.
In~\cite{Sirignano-Spiliopoulos:2018}, the solutions of a
high-dimensional PDE are approximated through a deep NN,
suitably trained to conform the behavior of a given differential
operator (equipped with initial and boundary conditions).
Physics-informed neural networks (PINN) are supervised learning
algorithms that constraint a numerical solution to satisfy a given set
of physics laws or PDEs, as studied for example in
\cite{Raissi-Perdikaris-Karniadakis:2018,%
  Kharazmi-Zhang-Karniadakis:2019,%
  Jin-Cai-Li-Karniadakis:2020,%
  Pang-DElia-Parks-Karniadakis:2020,%
  Shin-Darbon-Karniadakis:2020}.

%% MACHINE LEARNING, support vector machine 
As an alternative to the techniques mentioned above, we may  consider the
Support Vector Machine (SVM) approach.
Originally proposed at the beginning of the sixties
\cite{Vapnik-Lerner:1963,Vapnik-Chervonenkis:1964,Vapnik-Kotz:1982},
the SVM is a class of learning algorithms that can be used for solving
classification and regression problems, after a training over a
suitable set of input data~\cite{Vapnik:2000}.
The decision functions is a linear combination of a set of basis
functions that are nonlinear and parameterized by the so called
Support Vectors.
This approach can also provide a computationally efficient function
representation when the input space is high dimensional
~\cite{Guyon-Boser-Vapnik:1993,Suykens-Vandewalle:1999,Vapnik:1998,Scholkopf-Burges-Vapnik:1996}.
A SVM is trained to minimize the approximation errors, by searching for
a nonlinear function that fits the set of input data within an assigned
threshold and is as flat as possible.
The implementation is usually based on special 
polynomials or Gaussian kernels
\cite{Hofmann-Scholkopf-Smola:2008,Scholkopf-Smola:2009}.
An exhaustive review of SVM algorithms with their practical
applications is beyond the scope of this paper, so we refer the
interested reader to~\cite{Smola-Scholkopf:2004} and
\cite{Drucker-Burges-Kaufman-Smola-Vapnik:2003,CC01a} for a
more detailed introduction to both theory and usage.

The scope of this paper is to implement a couple of ML methods of the type
just described. As often happens in this kind of applications, one of
the most crucial issues will be the detection of the appropriate training sets.
The material is organized as follows.
%%
%% section 2
In Section~\ref{sec2:model:problem}, we briefly review the basic features
of Hermite spectral method for the 1D heat equation on the straight-line.
%%
%% section 3
In Section~\ref{sec3:quadrature} we discuss how the modification of
$\alpha$ may impact on the scheme formulation.
%%
%% section 4
In Section~\ref{sec4:spectral:algorithm} we reformulate the problem
as a collocation method, and introduce a discretization
in time admitting a scaling factor $\alpha=\alpha(t)$ changing with time.
%%
%% section 5
In Section~\ref{sec5:automatic:decision-making:approach} we design an
automatic decision-making strategy for the dynamical determination of $\alpha$,
applied to the numerical resolution of the heat equation with
homogeneous right-hand side. This is performed by training the ML
learning algorithm with the help of Gaussian-like profiles.
%%
%% section 6
In Section~\ref{sec6:nonhomogeneous:case} we apply a similar
strategy to the non-homogeneous case by employing a
training set containing suitable spline functions.
In Section~\ref{sec7:conclusions}, we provide some final remarks and discuss
future projects concerning the application of this ML strategy to the Vlasov equation.

% Hey Emacs, this is -*-latex-*-

%% sec2
\section{Galerkin approximation of the heat equation}
\label{sec2:model:problem}

Let $\REAL$ be the set of real numbers.
We work with the heat equation, that in strong form is given by:
\begin{subequations}
  \begin{align}
    \pt\us - \pxx\us &= \fs ,\phantom{\us_0(x)}  \qquad\xs\in\REAL,\quad\ts>0,\label{eq:strong:A}\\[0.5em]
    \us(x,0)         &= \us_0(x) ,\phantom{\fs}  \qquad\xs\in\REAL,           \label{eq:strong:B}  
  \end{align}
\end{subequations}
where $\fs$ is a given forcing term and $\us_0$ is the initial guess.
%% We set the diffusion parameter equal to one, since the choice of a
%% different value is not relevant to our exposition.
We assume that the solution $\us=\us (x,t)$ has an exponential decay at
infinity, although we do not specify at the moment the exact decay
rate.
This crucial issue will be discussed as we proceed with our
investigation.
For this reason, most of the passages that follow are to be intended
in informal way.

When $\fs=0$, the exact solution of
problem~\eqref{eq:strong:A}-\eqref{eq:strong:B} is explicitly known:
\begin{align}
  \us(x,t) = 2\alpha(t)w_\alpha (x,t)
  \quad\textrm{with}\ \ 
  \alpha(t) = \frac{1}{2\sqrt{t+1}}, \ \
  \wsa(x,t)=\exp\big(-(\alpha(t)\xs)^2\big).
  \label{eq:exact:solution}
\end{align}
As the reader can notice, the behavior of $\us$ for
$x\rightarrow\pm\infty$ depends on time.
This means that it is not easy to define an appropriate functional
space as a natural habitat for $\us$.
In particular, we would like to approximate $\us$
in~\eqref{eq:exact:solution} by Gaussian functions that display a
different decay rate (with $\alpha$ constant in time for instance).
For most of the paper we will continue to play with the simplified
case $\fs=0$, in order to prepare the ground for the non-homogeneous
case.

Let us go straight to the approximation.
For a positive integer number $\Ns$, let $\PS{\Ns}$ denote the space
of polynomials of degree at most $\Ns$.
Consider then the functional space:
\begin{align}
  \VsN(t) :=
  \Big\{
  \,\vsN(\cdot,t)=\wsa(\cdot,t)\phi_N\,|\,\phi_N\in\PS{N}\,
  \Big\},
  \label{eq:VsN:def}
\end{align}
where the weight function is given in~\eqref{eq:exact:solution}.
The classical Galerkin approximation is obtained in the usual way
through the following variational formulation:

\medskip
\textit{For $\ts\in[0,\Ts]$, find $\usN(t)\in\VsN (t)$ such that}:
\begin{align}
  &\int_{\REAL}\frac{\partial\usN}{\partial\ts}\phiN\dx
  +\int_{\REAL}\frac{\partial\usN}{\partial\xs}\,\frac{\partial\phiN}{\partial\xs}\dx
  =\int_{\REAL}\fs\phiN\dx ,
  \qquad\forall\phiN\in\PS{N},
  \label{eq:weak:form:IIa}
  \\[1em]
  &\int_{\REAL}\big(\usN(x,0)-\us_0(x)\big)\phiN(x,0)\dx=0 ,
  \qquad\forall\phiN\in\PS{N}.
  \label{eq:weak:form:IIb}
\end{align}

The exponential decrease of the weight function $\wsa(x,t)$ for
$\xs\to\pm\infty$ justifies the omission of the boundary terms in the
integration by parts.

Before proceeding, we need to introduce the set of Hermite
polynomials.
In fact, we will expand the approximation $\usN$ in the basis of
Hermite functions.
In the specific case, we first use Hermite polynomials multiplied by
the weight function $\wsa(\xs,\ts)$ introduced
in~~\eqref{eq:exact:solution}.
Later on, we will examine other options.

We recall that Hermite polynomials are recursively defined by
\begin{align}
  \Hs_{0} (\zeta) &= 1,\qquad\Hs_{1}(\zeta) = 2\zeta,\\[0.5em]
  \Hs_{n+1}(\zeta) &= 2\zeta\Hs_{n}(\zeta) - 2n\Hs_{n-1}(\zeta),\qquad\ns\geq2,
\end{align}
with the useful convention that $\Hs_{\ss}(\zeta)=0$ if $\ss<0$
~\cite{Shen-Tang-Wang:2011}.
These polynomials are orthogonal with respect to the weighted inner
product of $\LTWO(\REAL)$ using the weight function $\wsa$ with
$\alpha$ constantly equal to $1$.
In practice, we have:
\begin{align}
  \int_{\REAL}\Hs_{\ell}(\zeta)\Hs_{\ms}(\zeta)e^{-\zeta^2}\,d\zeta
  = 2^{\ms}\,\ms!\,\sqrt{\pi}\,\delta_{\ell,\ms},
  \label{eq:ortho:intg:00}
\end{align}
where $\delta_{\ell,\ms}$ is the Dirac delta symbol, whose values
equals one if $\ell=\ms$ and zero otherwise.
By the change of variable $\zeta=\alpha(t)\xs$, we are able to express
the orthogonality with respect to a general $\wsa$, which may depend
on time.
%%
%Further useful properties of Hermite polynomials are reported in the
% appendix in~\citeappx{in~\cite{Fatone-Funaro-Manzini:2021-arXiv}}.

From these assumptions, it is natural at this point to represent the
unknown $\usN$ through the expansion:
\begin{align}
  \usN(x,t) =
  \frac{\wsa(x,t)}{\sqrt{\pi}}\sum_{\ell=0}^{N}\ush_{\ell}(t)\Hs_{\ell}\big(\alpha(t)\xs\big),
  \label{eq:trial:function}
\end{align}
with its Fourier coefficients given by:
\begin{align}
  \ush_{m}(t) = \int_{\REAL}\usN(x,t)\phi_{m}(x,t)\,dx,
  \label{eq:Fourier:coeffs:def}
\end{align}
where
\begin{align}
  \phi_{m}(x,t) =
  \frac{\alpha(\ts)}{2^m\,m!}\Hs_{m}\big(\alpha(t)\xs\big) , \qquad
  0\leq\ms\leq\Ns.
  \label{eq:test:function}
\end{align}
In the trivial case where $f=0$, the series in~\eqref{eq:trial:function}
only contains the term corresponding to $\ell=0$, with $\ush_0=2\alpha (t)$
(see~\eqref{eq:exact:solution}).
In \eqref{eq:test:function}, the  choice for $\phi_{m}$, $0\leq\ms\leq\Ns$, is the one
also suggested for the test functions in~\eqref{eq:weak:form:IIa}.
Concerning the right-hand side, we have that the projection $\fs_N$ of
$\fs$ is obtained by:
\begin{align}
  \fsN(\xs,\ts)
  = \frac{\wsa(\xs,\ts)}{\sqrt{\pi}}\sum_{\ell=0}^{N}\fsh_{\ell}(\ts)\Hs_{\ell}(\alpha(\ts)\xs).
\end{align}
The Fourier coefficients $\fsh_{\ms}(\ts)$ are recovered after
multiplying by the test functions in~\eqref{eq:test:function}
and integrating over $\REAL$, so obtaining:
%%
%% ELSEVIER
\begin{align}
  &\frac{\alpha(\ts)}{2^{\ms}\,\ms!}\int_{\REAL}\fsN(\xs,\ts)\Hs_{\ms}(\alpha(\ts)\xs),\dx
  = \frac{1}{2^m\,\ms!\,\sqrt{\pi}}
  \sum_{\ell=0}^{N}\fsh_{\ell}(\ts)
  \int_{\REAL}\wsa(\xs,\ts)\Hs_{\ell}(\alpha(\ts)\xs)\Hs_{\ms}(\alpha(\ts),\xs),\alpha(\ts)\dx
  %% ---------------
  \nonumber\\[0.5em]
  &= \frac{1}{2^m\,\ms!\,\sqrt{\pi}}
  \sum_{\ell=0}^{N}\fsh_{\ell}(\ts)
  \int_{\REAL}\exp{(-\zeta^2)}\Hs_{\ell}(\zeta)\Hs_{\ms}(\zeta)\,d\zeta
  %% ---------------
  %% \nonumber\\[0.5em]
  %% &
  = \frac{1}{2^m\,\ms!\,\sqrt{\pi}}
  \sum_{\ell=0}^{N}\fsh_{\ell}(\ts)\,2^{\ms}\ms!\sqrt{\pi}\delta_{\ell,\ms}
  %% ---------------
  %% \nonumber\\[0.5em]
  %% &
  = \fsh_{\ms}(\ts),
  \label{eq:RHS:function}
\end{align}
where we used the orthogonality properties of Hermite polynomials.

We can go through the computations in order to find out the scheme in
finite dimension.
The procedure is here omitted, though it follows from straightforward
(but quite annoying) calculations.
Some passages are briefly reported in the 
appendix~\citeappx{in~\cite{Fatone-Funaro-Manzini:2021-arXiv}}.
We then arrive at the following result.

{\it Let $\ush_{\ms}(\ts)$ and $\fsh_{\ms}(\ts)$ be the $\ms$-th
time-dependent Fourier coefficients of $\usN$ and $\fsN$, respectively
introduced in~\eqref{eq:trial:function} and~\eqref{eq:RHS:function}.
Then, the Fourier coefficients $\ush_{\ms}(\ts)$, $0\leq\ms\leq\Ns$,
satisfy the system of ordinary differential equations:}
\begin{align}
  \pt\ush_{\ms}(\ts) -
  \bigg(\frac{\asp(\ts)}{2\alpha(\ts)} + \alpha(\ts)^2\bigg)\ush_{\ms-2}(\ts)
  -(\ms+1)\frac{\asp(\ts)}{\alpha(\ts)}\ush_{\ms}(\ts)
  = \fsh_{\ms}(\ts) .
  \label{eq:form0:final:15}
\end{align}

The same relation between the coefficient of the Fourier expansions is
obtained from the following variational formulation of the heat
equation:
\begin{align}
  &\int_{\REAL}
  \Big(\frac{\partial\us}{\partial\ts}-\frac{\partial^2\us}{\partial\xs^2}\Big)\phi\dx
  = \int_{\REAL}\fs\phi\dx
  \quad\forall\phi\,(\mbox{generic test function}),
  \label{eq:weak:form:Ia}
\end{align}
where, for the moment, $\phi$ denotes a generic test function, and no
integration by parts has been performed.
The equivalence of the two versions is formally shown in the
appendix~\citeappx{in~\cite{Fatone-Funaro-Manzini:2021-arXiv}.}
Nevertheless, the last formulation cannot be easily embedded in the
proper functional spaces, so that an analysis of convergence for the
corresponding Galerkin approximation is hard to achieve.
A third version has been proposed in~\cite{Ma-Sun-Tang:2005}.
This is based on a rewriting of the heat equation as follows:
\begin{align}
  &\frac{d}{\dt}\int_{\REAL}\us\phi\dx
  -\int_{\REAL}\us\,\Big(\frac{\partial\phi}{\partial\ts}
  +\frac{\partial^2\phi}{\partial\xs^2}\Big)\dx
  =\int_{\REAL}\fs\phi\dx
  \qquad\forall\phi\,(\mbox{generic test function}).
  \label{eq:cinesi}
\end{align}

\noindent
%% marco-comment
%\textbf{Bisogna ricontrollare quello che diciamo nel paragrafo che
%  segue!}
%%
This leads again to the scheme~\eqref{eq:form0:final:15}.
The reason for this last choice is due to an erroneous writing of the
relation linking the Fourier coefficients of the classical Galerkin
scheme \eqref{eq:weak:form:IIa}, which led
to equation (1.4) in~\cite{Ma-Sun-Tang:2005}.
This flaw suggested to the authors a modification of the variational
formulation according to~\eqref{eq:cinesi}.
After a review of the computations, we consider the alternative
version~\eqref{eq:cinesi} unnecessary.
Details of the  passages are provided in the 
appendix~\citeappx{in~\cite{Fatone-Funaro-Manzini:2021-arXiv}}.
It has to be remarked however that the work in~\cite{Ma-Sun-Tang:2005}
has been mainly focused to the study of convergence of the
scheme~\eqref{eq:form0:final:15}, which is exactly the one we are
implementing in the present paper.
Therefore, the theoretical results developed
in~\cite{Ma-Sun-Tang:2005} are definitely interesting to us.

The exact solution \eqref{eq:exact:solution} also satisfies the
recursive relation in~\eqref{eq:form0:final:15} for $\fs=0$.
In fact, we first recall that $\alpha(t)=1\slash{2\sqrt{t+1}}$, which
implies $\asp(t)\slash{2\alpha(t)}=-\alpha(t)^2$.
In this way,~\eqref{eq:form0:final:15} reduces to:
\begin{align}
  \pt\ush_{\ms}(\ts) + 2(\ms+1)\alpha(\ts)^2\ush_{\ms}(\ts) = \fsh_{\ms}(\ts).
  \label{eq:form0:final:150}
\end{align}
All the Hermite coefficients of $\us$ are zero with the exception of
the one corresponding to $m=0$, where we have
$\hat\us_0(t)=2\alpha(t)$.
Relation~\eqref{eq:form0:final:150} is actually compatible with this
choice.

Let us now suppose that the parameter $\alpha$ is constant.
In this circumstance, always for $\fs=0$, we can find for example an
explicit expression for the coefficients in~\eqref{eq:form0:final:15}
for $m\geq 0$ ($m$ even):
\begin{align}
  \ush_{m}(t) &= \frac{\sqrt{\pi}\alpha^m}{\left(\frac{m}{2}\right)!}t^{m/2},
  \label{eq:Fourier:coeff}
\end{align}
from which we get $\ush_0=\sqrt{\pi}$, which corresponds to the
initial datum $\us(\xs,0)=\wsa(\xs,0)$.
Indeed, by setting $\alpha$ constant,
relation~\eqref{eq:form0:final:15} becomes:
\begin{align}
  \pt\ush_{\ms}(\ts) - \alpha^2\ush_{\ms-2}(\ts) = 0.
  \label{eq:ode:constant:alpha}
\end{align}
Such an equation is verified by noting that:
\begin{align*}
  \pt\ush_{m}(t)
  = \frac{\sqrt{\pi}\alpha^m}{\left(\frac{m}{2}\right)!}\pt t^{m/2}
  %% = \frac{\sqrt{\pi}\alpha^m}{\left(\frac{m}{2}\right)!}\frac{m}{2}t^{\frac{m-2}{2}}
  %% = \frac{\sqrt{\pi}\alpha^{m-2}}{\left(\frac{m}{2}-1\right)!}\alpha^2t^{({m-2})/{2}}
  = \alpha^2\frac{\sqrt{\pi}\alpha^{m-2}}{\left(\frac{m-2}{2}\right)!}
  t^{({m-2})/{2}}
  = \alpha^2\ush_{m-2}.
\end{align*}
%We would like to shqow that the function constructed from the
%coefficients in~\eqref{eq:Fourier:coeff} does not satisfy the heat
%equation~\eqref{eq:strong:A} for $\alpha=1/2$.
Since, for the initial datum in~\eqref{eq:strong:B}, only the first
mode is different from zero, we have that $\us_0$ is an even function
that remains even during all the time evolution (recall that $\fs=0$).
Therefore, by using the Fourier coefficients
in~\eqref{eq:Fourier:coeff}, we build the sum: 
%% marco-comment
%\textbf{Attenzione, da
%  qui in poi scriviamo $\wsa(\xs)$ invece che $\wsa(\xs,\ts)$}{\color{red}perche non dipende da t}
%%
\begin{align}
  \vs(\xs,\ts)
  = w_\alpha(x)
 % \exp\big(-(\alpha\xs)^2\big)
  \sum_{m(even)=0}^{\infty}\alpha^{m}\frac{t^{\frac{m}{2}}}{(\frac{m}{2})!}\Hs_{m}(\alpha\xs)= w_\alpha(x)
  %\exp\big(-(\alpha\xs)^2\big)
  \sum_{\ell =0}^{\infty}\alpha^{2\ell}\frac{t^\ell}{\ell!}\Hs_{2\ell}(\alpha\xs),
  \label{eq:exact:solution:constant:alpha}
\end{align}
where we dropped the dependence on $t$ in $\wsa$, since the weight no longer depends on time.
We have that $\vs(\xs,0)=\us(\xs,0)$ when $\alpha =1/2$.

We can also say that $v(0,t)=u(0,t)$, for $0\leq t<1$.
We prove this fact by writing the MacLaurin expansion:
\begin{align}
  \us(0,t)
  = 2\alpha(t)
  = \frac{1}{\sqrt{t+1}}
  = \sum_{\ell=0}^{\infty} \frac{(-1)^\ell}{2^\ell\ell!}\Big(\prod_{k=1}^{\ell}(2i-1)\Big)t^{\ell}.
  \label{eq:maclaurin:expansion:alpha}
\end{align}
Now, the expansions \eqref{eq:exact:solution:constant:alpha} and
\eqref{eq:maclaurin:expansion:alpha} are the same, as it can be
checked by expressing the value of $\Hs_{2\ell}(0)$ through the
recursion formula for Hermite polynomials.
Let us observe that the above mentioned series have a convergence
radius equal to 1, therefore they diverge for $t> 1$.
%%
%However, if we assume the right dependence on time of the coefficient
%$\alpha(t)$, which is the one given in~\eqref{eq:exact:solution}, we
%can use the scheme in Proposition~\ref{prop:numerical:scheme} to
%approximate solution \eqref{eq:exact:solution}, as is done
%in~\cite{CINESI}.

We are tempted to deduce that the functions $\vs$ and $\us$ always
coincide both in time and space.
This kind of equality is not clear however, since it presumes a
correct interpretation of the type of convergence of the series
in~\eqref{eq:exact:solution:constant:alpha}.
For any fixed $\xs$ and $\ts$, we may expect point-wise convergence,
but it is hard to conclude that this is going to be true in some
normed functional space.
In fact, we are trying to represents a Hermite function with a certain
Gaussian decay, through an expansion by other Hermite functions
associated to a different Gaussian behavior.
To this respect, let us note that, by choosing $\alpha(t)$ as in
\eqref{eq:exact:solution}, only the first mode ($m=0$) is activated
during the whole computation, so that the numerical error exclusively
depends on the time discretization procedure.
On the contrary, for $\alpha$ fixed equal to $1/2$, all the Fourier
coefficients are involved.

In~\cite{Ma-Sun-Tang:2005}, spectral convergence is proven for a
special choice of $\alpha (t)$ (see relation (2.1) in that paper),
which includes the case \eqref{eq:exact:solution}.
The analysis for a more general $\alpha$ is at the moment unavailable.
Such an extension is rather troublesome and it is due to the
difficulty to find proper Sobolev type inclusions between spaces
displaying different decay behaviors at infinity of the weight
functions.
The results of the experiments obtained by truncating the sum
in~\eqref{eq:exact:solution:constant:alpha} at a fixed integer $N$
show an excellent agreement with the exact solution $u$, for $t<1$.
As $t$ reaches the value 1, oscillations are observed. Like in the
Gibbs' phenomenon, these do not diminish by increasing $N$
(recall that the convergence radius of the involved series is equal to 1).

%% ----
%% sec3
%%-----

\section{Quadrature nodes and other useful formulas}
\label{sec3:quadrature}

The considerations made in the previous section, rise the question of
the passage from the Fourier coefficients relative to a certain choice
of $\alpha$ to those related to a different value of such a parameter.
From the practical viewpoint, we briefly provide a recipe for this
kind of basis transformation.
%Now, we consider the problem of how to change the weighting
%coefficient from a given value of $\alpha$, e.g., the value taken by
%at a given time $t$, to another value, for example $1$.
%%
Suppose that we are given a function $\vs$ expressible as a series of
Hermite functions for some $\alpha>0$:
\begin{align}
  \vs(\xs) =
  \frac{ w_{\alpha}(\xs) }{\sqrt{\pi}}
  \sum_{\ell=0}^{\infty}\vsh_{\ell}
  \Hs_{\ell}\big(\alpha\xs\big),
  \label{eq:vs:alpha:def}
\end{align}
with Fourier coefficients
\begin{align}
  \vsh_{m} = \frac{\alpha}{2^m\,m!}
  \int_{\REAL}
  \vs(\xs)
  \Hs_{m}(\alpha\xs)
  \,dx=
   \frac{1}{2^m\,m!}
  \int_{\REAL}
  \vs(x/\alpha)
  \Hs_{m}(\xs)
  \,dx ,
  \qquad m\geq 0.
  \label{eq:vs:alpha:Fourier:coeffs:def}
\end{align}
We would like to express the same function as the series of Hermite
functions relative to $\alpha =1$:
\begin{align}
  \vs(\xs) =
  \frac{ w_1(\xs) }{\sqrt{\pi}}
  \sum_{\ell=0}^{\infty}\Vsh_{\ell}
  \Hs_{\ell}(\xs),
  \label{eq:vs:unit:def}
\end{align}
with Fourier coefficients
\begin{align}
  \Vsh_{m} = \frac{1}{2^m\,m!}\int_{\REAL}
  \vs(\xs)
  \Hs_{m}(\xs)
  \,dx ,
  \qquad  m\geq 0.
  \label{eq:vs:unit:Fourier:coeffs:def}
\end{align}
We can immediately get a formula that provides the set of Fourier
coefficient $\Vsh_{m}$ from the other one.
It is enough to informally substitute
expansion~\eqref{eq:vs:alpha:def} in the last integral.
So, we find that
\begin{align}
  \Vsh_{m} = \frac{1}{2^m\,m!\sqrt{\pi}}
  \sum_{\ell=0}^{\infty}
  \vsh_{\ell}
  \int_{\REAL}
  \Hs_{\ell}\big(\alpha\xs\big)
  \Hs_{m}(\xs)
  w_\alpha(\xs)
  %\exp\big(-(\alpha\xs)^2\big)
  \,dx =
  \frac{1}{\alpha 2^m\,m!\sqrt{\pi}}
  \sum_{\ell=0}^{\infty}
  \vsh_{\ell}
  \int_{\REAL}
  \Hs_{\ell}\big(\xs\big)
  \Hs_{m}\left(\frac{\xs}{\alpha}\right)
  w_1(x)
  \,dx
  \label{eq:vsh:to:Vsh}
\end{align}
and, viceversa, substituting~\eqref{eq:vs:unit:def}
in~\eqref{eq:vs:alpha:Fourier:coeffs:def}.
\begin{align}
  \vsh_{m} = \frac{\alpha}{2^m\,m!\sqrt{\pi}}
  \sum_{\ell=0}^{\infty}
  \Vsh_{\ell}
  \int_{\REAL}
  \Hs_{\ell}(\xs)
  \Hs_{m}(\alpha\xs)
  w_1(\xs)
  %  \exp\big(-\xs^2\big)
  \,dx ,
  \qquad m\geq 0.
  \label{eq:Vsh:to:vs}
\end{align}

Explicit formulas are available for the computation of these
integrals.
The procedure looks however complex and expensive, so that we provide
here below an alternative by using numerical quadrature.
As usual, for a fixed integer $\Ns$, we introduce the zeros of
$\Hs_N$.
These nodes will be denoted by $\xiN_j$, $j=1,\ldots,\Ns$.
Together with the weights~\cite[Eq.~(3.4.9)]{Funaro:1992}:
\begin{align}
  \wsN_j = \sqrt{\pi}2^{N+1}N!\big[\Hsp_{n}(\xiN_j)\big]^{-2} ,
  \qquad
  1\leq\js\leq\Ns,
  \label{eq:pesi:quadratura}
\end{align}
%%
%% formula di quadratura
we have the quadrature formula:
\begin{align}
  \int_{\REAL}\ps(\xs)dx = \sum_{j=1}^N\ps(\xiN_j)\wsN_j,
  \label{eq:quadratura}
\end{align}
which is exact for any polynomial $\ps$ of degree less than $2N$.
We can pass from the point-values to the Fourier coefficients by
observing that ~\eqref{eq:vs:alpha:Fourier:coeffs:def} can be
rewritten as:
\begin{align}
  \vsh_{m} = 
  \frac{1}{2^m\,m!}
  \sum_{j=1}^{\Ns}v(\xiN_j/\alpha)H_m(\xiN_j)\wsN_j,
  \label{eq:vs:alpha:Fourier:coeffs2:def}
\end{align}
provided the product $\vs\Hs_m$ is a polynomial of degree less than
$2N$.
A direct application of the quadrature rule to~\eqref{eq:vsh:to:Vsh}
yields:
\begin{align}
  \Vsh_{m}
  \approx
  \frac{1}{\alpha 2^m\,m!\sqrt{\pi}}
  \sum_{\ell=0}^{\Ns-1}
  \vsh_{\ell}
  \sum_{j=1}^{N}
  \Hs_{\ell}\big( \xiN_j\big)
  \Hs_{m}\left({\xiN_j}/{\alpha}\right)\wsN_j ,
  \qquad 0\leq\ms\leq N-1,
  \label{eq:vsh:to:Vsh:quadrature}
\end{align}
where the sum on $\ell$ has been suitably truncated.
The change of coefficients from a basis to another one is then
obtained at a reasonable cost, although some approximation has been
introduced.
The writing in~\eqref{eq:vsh:to:Vsh:quadrature} may be represented by
a linear operator associated with an $\Ns\times\Ns$ matrix.
It must be observed that, once $N$ has been fixed, we can get reliable
approximations only if $\alpha$ stays within a certain range.
In fact, it is necessary to have an appropriate amount of nodes in the
support of the Gaussian (we practically define the ``support" as the set where the Gaussian
is different from zero for a given machine error precision). 
If $\alpha$ is too large, the support is narrow and all the nodes are
outside the region of interest.
In this case, the interpolating polynomial produces oscillations.
On the contrary, for $\alpha$ too small, the support is wide and the
nodes turn out to be concentrated at its interior, so providing a poor
approximation.
We will return to these issues in
Section~\ref{sec6:nonhomogeneous:case}.

The use of quadrature formulas suggests the collocation method as a
possible alternative to the Galerkin scheme~\eqref{eq:form0:final:15}.
In order to follow this new path, it is customary to introduce the
polynomial Lagrange basis:
\begin{align*}
  \lsN_j(\xs) =
    \frac{\Hs_{\Ns}(\xs)}{\Hsp_{\Ns}(\xiN_j)(\xs-\xiN_j)}, \ \ {\rm if}~\xs\neq\xiN_j,
   \qquad\qquad \lsN_j(\xiN_j) = 1   ,                                          \end{align*}
so that, for any polynomial $\ps$ of degree at most $N-1$ one has:
\begin{align*}
  \ps(\xs) = \sum_{j=1}^{\Ns}\ps(\xiN_j)\lsN_j(\xs).
\end{align*}
We also recall the $\Ns\times\Ns$ differentiation matrices
$\Ds^{(1)}=\big\{\ds^{(1)}_{ij}\big\}$ and
$\Ds^{(2)}=\big\{\ds^{(2)}_{ij}\big\}$, with
$\Ds^{(2)}=\Ds^{(1)}\Ds^{(1)}$.
These are recovered from the first and second derivatives of the
Lagrange polynomials:
\begin{align}
  \ds^{(1)}_{ij} =
  \frac{d\lsN_j}{\dx}(\xiN_i) ,\qquad\qquad
  \ds^{(2)}_{ij} =
  \frac{d^2\lsN_j}{\dx^2}(\xiN_i).
  \label{eq:vsh:to:Vsh:derivate}
\end{align}
The explicit entries of these matrices are well known and can be found
in many textbooks, e.g.,~\cite{Funaro:1992}
and~\cite{Shen-Tang-Wang:2011}.
%% {\color{red} mio libro [7]} or {\color{red} libro Shen Tao Wang}.
We now have all the elements to construct an approximation method of
collocation type.
This will be developed in the next section.

%% ----
%% sec4
%% ----

%% newcommands da spostare nel file defs.tex
\newcommand{\wgh}{\omega} %% Hermite'sweight

\section{Collocation algorithm and time discretization}
\label{sec4:spectral:algorithm}

We start by writing the unknown $u$ as:
\begin{align}
  \us(\xs,\ts)
  = \ps(\xs,\ts)\wgh_{\alpha}(\xs,\ts),
  \label{eq:up}
\end{align}
where we recall that
$\wgh_{\alpha}(\xs,\ts)=\exp\big(-\alpha(t)^2\xs^2\big)$.
Note that in~\eqref{eq:up}  $\ps$ is not necessarily a polynomial.
By substituting this expression into the heat equation we obtain:
\begin{align}
  \wgh_{\alpha}\pt\ps + \ps\pt\wgh_{\alpha}
  = \wgh_{\alpha}\pxx\ps - 2\px\ps\,\px\wgh_{\alpha} + \ps\pxx\wgh_{\alpha} + \fs.
  \label{eq:passaggio1}
\end{align}
We then divide both sides of~\eqref{eq:passaggio1}  by $\wgh_{\alpha}$ to obtain:
\begin{align}
  \pt\ps + 2\alpha\alpha'\xs^2\ps
  = \pxx\ps - 4\alpha^2\xs\px\ps
  + 2\alpha^2\big(2\alpha^2\xs^2-1\big)\ps
  + \frac{\fs}{\wgh_{\alpha}},
  \label{eq:passaggio2}
\end{align}
where $\alpha'$ denotes the derivative of $\alpha$ with respect to
time.
Finally, we look for an approximating polynomial $p_N$ of degree less
than or equal to $N-1$, after collocating the above equation at the
points $\xiN_j$, $1\leq j\leq N$:
\begin{multline}
  \pt\ps_N(\xiN_j,\ts) + 2\alpha(\ts)\alpha'(\ts)(\xiN_j)^2\ps_N(\xiN_j,\ts)
  = (\pxx\ps_N)(\xiN_j,\ts)
  \\[0.25em]  %{\color{red}accorciare \ sopra} 
  - 4\alpha^2(\ts)\xiN_j(\px\ps_N)(\xiN_j,\ts)
  + 2\alpha^2(\ts)\big(2(\alpha (\ts)\xiN_j)^2 - 1\big)\ps_N(\xiN_j,\ts)
  + \frac{\fs(\xiN_j,\ts)}{\wgh_{\alpha}(\xiN_j,\ts)}.
  \label{eq:collocazione}
\end{multline}
Here, the derivatives $\px\ps_N$ and $\pxx\ps_N$ can be evaluated with
the help of the finite dimensional operators $D^{(1)}$ and $D^{(2)}$
defined in ~\eqref{eq:vsh:to:Vsh:derivate}.
Note that the approximate solution of the heat equation is recovered
as in ~\eqref{eq:up} through the expression
\begin{equation}
  u_N(x,t)=p_N(x,t)w_{\alpha}(x,t).
\end{equation}
If $\alpha (t)$ is chosen to follow the behavior of the exact solution
$u$ for $f=0$, i.e. $\alpha (t)=1/2\sqrt{t+1}$ as
in~\eqref{eq:exact:solution}, the approximation is excellent.
Indeed, in this simple case, $p_N$ turns out to be a polynomial of
degree zero in the variable $x$, i.e., $p_N(x,t)=\alpha (t)$, as it can be
checked by direct substitution.

Since the exact $\alpha (t)$ is not known a priori, it is interesting
to see what happens for different choices of this parameter.
As in the Galerkin method, we start by considering the case of
$\alpha$ constant in time.
Thus, the scheme \eqref{eq:collocazione} becomes:
%%
%% SIAM
%%
\begin{align}
  \pt\ps_N(\xiN_j,\ts) 
  = (\pxx\ps_N)(\xiN_j,\ts)
  - 4\alpha^2\xiN_j(\px\ps_N)(\xiN_j,\ts)
  + 2\alpha^2\big(2(\alpha\xiN_j)^2 - 1\big)\ps_N(\xiN_j,\ts)
  + \frac{\fs(\xiN_j,\ts)}{\wgh_{\alpha}(\xiN_j)}.
  \label{eq:collocazione2}
\end{align}
The next step is to discretize in time for $\ts\in[0,\Ts]$.
To this end we just apply the forward Euler method, though we are
conscious that more sophisticated methods are available.
We denote the time step by $\Delta t$  and the generic $n$-th point of the
time grid by $\ts^n=n\Delta t$, $0\le n \le {\cal N}$, $t^0=0$, ${\cal N} \Delta t=T$.
We obtain
%%
%% ELSEVIER
\begin{align}
  \label{eq:eulerocollocazione}
  \frac{\ps_N^{n+1}(\xiN_j) -\ps_N^{n}(\xiN_j)}{\Delta t} 
  = (\pxx\ps_N^n)(\xiN_j)
  - 4\alpha^2\xiN_j(\px\ps_N^n)(\xiN_j)
  + 2\alpha^2\big(2(\alpha\xiN_j)^2 - 1\big)\ps_N^n(\xiN_j)
   + \frac{\fs(\xiN_j,\ts^n)}{\wgh_{\alpha}(\xiN_j)}.
\end{align}
We supplement equation \eqref{eq:eulerocollocazione} with the initial
condition $\ps_N^0$, which is derived from the initial solution
$\us(\xs,0)$.

Again we take for the moment $f=0$, in order to perform some
experiments.
We would like to have the chance to modify $\alpha$ during the
evolution.
The idea is to keep it constant for a certain number of iterations,
and then shift to another constant value.
The problem is that the piecewise constant function $\alpha (t)$ is
not differentiable at the discontinuity points,
%\textbf{non mi piace
%  ``\emph{discontinuity}''}{\color{red}il latex lo ritiene corretto, e anche l'analista}, 
so that the scheme
\eqref{eq:collocazione} cannot be applied, as it requires the
knowledge of $\alpha'$.
We then modify the time advancing procedure in the following fashion.

%in every time
%interval $\big[\ts^n,\ts^{n+1}\big]$.
%%
Starting from a value $t^n$ of the time discretization grid, we keep
$\alpha$ constantly equal to a given $\alpha_n$, during $m$ consecutive time steps,
i.e. within the interval $\big[\ts^n,\ts^{n+m}\big]$.
In this way, we find $p^n_N, \cdots p^{n+m}_N$ solving equation
\eqref{eq:eulerocollocazione}.
Correspondingly, we also construct:
$u_N^{n+k}=p_N^{n+ k}\exp\big(-\alpha_{n}^2x^2\big)$,  $0\le  k \le m$.
At time $\ts^{n+m}$ we decide to assign a new $\alpha=\alpha_{n+m}$ (we will
see later how this can be done automatically).
Before continuing with the scheme \eqref{eq:eulerocollocazione}, we
need to pass from the representation of $u^{n+m}_N$, with the Gaussian
associated to $\alpha_n$, to the new one relative to $\alpha_{n+m}$.
In order to do this, from $p^{n+m}_N$ we build a new polynomial such
that:
\begin{align}
  \widetilde\ps_N(\xiN_j) = \ps^{n+m}_N(\xiN_j)\exp\big( ( \alpha_{n+m}^2 - \alpha_{n}^2 )(\xiN_j)^2 \big),
  \qquad 1\leq j\leq N.
  \label{eq:interpolazione}
\end{align}
This transitory $\widetilde p_N$ is used as a new initial guess to advance
in time, until the next update of $\alpha$. Of course, if
$\alpha_{n+m}=\alpha_n$, this last passage is not relevant.

We now discuss the results of a series of numerical experiments, where
the exact solution $u$ is compared with its approximation $u_N$.
The errors are evaluated in the followings norms:
%%
%% ELSEVIER
\begin{align}
  {\cal N}^1=\left( \sum_{j=0}^N\vert u(\xiN_j)-u_N(\xiN_j)\vert^2\right)^{\hspace{-.2cm}1/2},\qquad
  {\cal N}^2=\max_\REAL\vert u-u_N\vert ,\qquad
  {\cal N}^3=\left( \sum_{j=0}^N\vert u(\xiN_j)-u_N(\xiN_j)\vert^2 \wsN_j\right)^{\hspace{-.2cm}1/2}
  \label{eq:norme}
\end{align}
with ${\cal E}_j=\vert u(\xiN_j)-u_N(\xiN_j)\vert$, $1\leq j\leq N$.
The maximum in ${\cal N}^2$ is computed on a very fine grid on the
$x$-axis and $\wsN_j$, $1\leq j\leq N$ are the quadrature weights
introduced in~\eqref{eq:pesi:quadratura}.
The time step has been chosen small enough in order not to affect the
error in space.
Since the sum of the weights is equal to $\sqrt{\pi}$, the error
evaluated using norm ${\cal N}^3$ is less than the error evaluated
using norm ${\cal N}^1$.
Therefore, with the exception of Table \ref{table:T1}, we will not
report the results related to the last norm in~\eqref{eq:norme}.

\begin{figure}[ht!]
  \centering
  \includegraphics[width=0.7\linewidth]{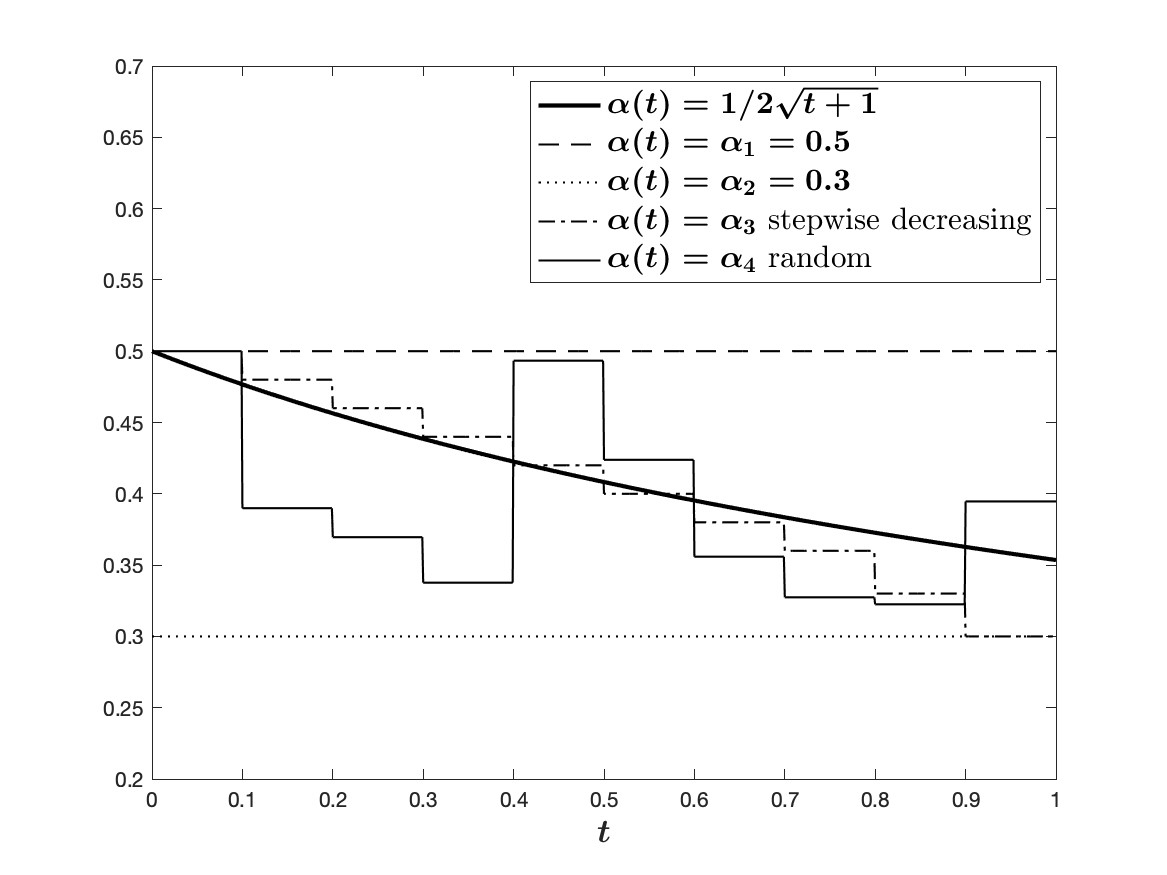}%%{./Figures/figAlpha1.jpg}
  \vskip-0.4truecm
  \caption{ Choices of $\alpha$ in the various experiments.}
  \label{fig:a_vero} 
\end{figure}

We continue to consider $f=0$ and we integrate in time for
$\ts\in[0,\Ts]$ with $\Ts=1$, $\Delta t= 10^{-7}$.
The first tests in Table \ref{table:T1}, show a spectral decay of the
solution when $\alpha =\alpha (t)$ is taken as in
\eqref{eq:exact:solution}.
Since we do not know a priori the behavior of $\alpha$,
the interesting part comes when we choose this parameter according to some
prescribed law.
In Table \ref{table:T2} we find the results where $\alpha$ has been
fixed once and for all.
As expected, the errors are not encouraging.  An improvement is
obtained by imposing that $\alpha (t)$ is piecewise constant.
To this end, we subdivide the time interval
%% $[0,1]$
in ten subintervals.
In Table \ref{table:T3} (left panel) we have the expression of the
errors when $\alpha$ starts from the value $0.5$ and decreases up to
$0.3$ in 9 equal steps.
In Table \ref{table:T3} (right panel), ten guesses of $\alpha$ have
been made randomly in the interval $[0.3, 0.5]$.
The behavior of $\alpha$ in all the above examples is summarized in
Fig.~\ref{fig:a_vero}.

%% ELSEVIER
%% \input{sec4_tables_ELSEVIER.tex}
\renewcommand{\TABROW}[3]{ \hspace{0.3cm} #1 \hspace{0.4cm} & \hspace{0.4cm} #2 \hspace{0.4cm} & \hspace{0.4cm} #3 \hspace{0.65cm}}

%%%%%%%%% alpha=alpha Vero %%%%%%%%%

\renewcommand{\TABROW}[4]{ #1 & #2 & #3 & #4}
\begin{table}[ht!]
  \centering
  \begin{tabular}{l|c|c|c}
    \hline\hline
    \TABROW{$N$}{$\displaystyle\calN^1$}{$\displaystyle\calN^2$}{$\displaystyle\calN^3$}\\%[0.25em]
    \hline\hline
    \TABROW{ 4}{2.6171e-04}{2.2846e-04}{1.0907e-04}\\
    \TABROW{ 6}{1.2092e-05}{1.2562e-05}{2.2202e-06}\\
    \TABROW{ 8}{6.2025e-07}{7.0862e-07}{5.5019e-08}\\
    \TABROW{10}{3.0252e-08}{4.1672e-08}{7.6304e-09}\\
    \hline\hline                  
  \end{tabular}
  \caption{Errors at time $T=1$ (using $\Delta t= 10^{-7}$) and $\alpha
    (t)=1/2\sqrt{t+1}$, $t\in [0,1]$.}
  \label{table:T1}
\end{table}
\renewcommand{\TABROW}[3]{ \hspace{0.3cm} #1 & \hspace{0.3cm} #2 & \hspace{0.3cm} #3}
\begin{table}[ht!]
  \begin{minipage}[b]{0.4\linewidth}
    \hspace{3cm}
    %%%%%%%%% alpha=0.5 %%%%%%%%%
    \begin{tabular}{l|c|c}
      \hline\hline
      \TABROW{$N$}{$\displaystyle{\cal N}^1$}{$\displaystyle{\cal N}^2$}\\%[0.5em]
      \hline\hline
      \TABROW{ 4}{ 2.6090e-02 }{ 6.8507e-02 } \\[-0.5em]
      \TABROW{ 6}{ 8.7970e-03 }{ 2.6723e-03 } \\[-0.5em]
      \TABROW{ 8}{ 3.0084e-03 }{ 1.0619e-02 } \\[-0.5em]
      \TABROW{10}{ 1.0421e-03 }{ 4.2718e-03 } \\%%[-0.5em]
      \hline\hline                  
    \end{tabular}
  \end{minipage}
  %% \hspace{1cm}
  \begin{minipage}[b]{0.4\linewidth}
    \hspace{2cm}
    %%%%%%%%% alpha=0.3 %%%%%%%%%
    \begin{tabular}{l|c|c}
      \hline\hline
      \TABROW{$N$}{$\displaystyle{\cal N}^1$}{$\displaystyle{\cal N}^2$}\\%[0.5em]
      \hline\hline
      \TABROW{ 4}{ 6.6306e-02 }{ 8.0850e-02 } \\[-0.5em]
      \TABROW{ 6}{ 4.8930e-02 }{ 7.2890e-02 } \\[-0.5em]
      \TABROW{ 8}{ 4.0269e-02 }{ 7.5228e-02 } \\[-0.5em]
      \TABROW{10}{ 3.5001e-02 }{ 8.3357e-02 } \\%%[-0.5em]
      \hline\hline                  
    \end{tabular}
  \end{minipage}
  \vspace{2mm}
  \caption{Errors at time $T=1$ (using $\Delta t= 10^{-7}$) and $\alpha$
    constant in time: $\alpha =0.5$ (left); $\alpha =0.3$ (right).}
  \label{table:T2}
\end{table}

%%%%%%%%% alpha stepwise %%%%%%%%%
\begin{table}[ht!]
  \begin{minipage}[b]{0.4\linewidth}
    \hspace{3cm}
    %%%%%%%%% alpha=0.5 %%%%%%%%%
    \begin{tabular}{l|c|c}
      \hline\hline
      \TABROW{$N$}{$\displaystyle{\cal N}^1$}{$\displaystyle{\cal N}^2$}\\%[0.5em]
      \hline\hline
      \TABROW{ 4}{ 1.4775e-03 }{ 2.0235e-02 } \\[-0.5em]
      \TABROW{ 6}{ 8.6681e-05 }{ 5.2286e-03 } \\[-0.5em]
      \TABROW{ 8}{ 1.6311e-05 }{ 1.5659e-03 } \\[-0.5em]
      \TABROW{10}{ 1.3840e-06 }{ 4.6910e-04 } \\%%[-0.5em]
      \hline\hline                  
    \end{tabular}
  \end{minipage}
  %%\hfill
  \begin{minipage}[b]{0.4\linewidth}
    \hspace{2cm}
    %%%%%%%%% alpha=0.3 %%%%%%%%%
    \begin{tabular}{l|c|c}
      \hline\hline
      \TABROW{$N$}{$\displaystyle{\cal N}^1$}{$\displaystyle{\cal N}^2$}\\%[0.5em]
      \hline\hline
      \TABROW{ 4}{ 9.9108e-03 }{ 1.2908e-02 } \\[-0.5em]
      \TABROW{ 6}{ 2.6628e-03 }{ 2.2971e-03 } \\[-0.5em]
      \TABROW{ 8}{ 9.1279e-04 }{ 1.1427e-03 } \\[-0.5em]
      \TABROW{10}{ 3.0883e-04 }{ 3.9922e-04 } \\%%[-0.5em]
      \hline\hline                  
    \end{tabular}
  \end{minipage}
  \vspace{2mm}
  \caption{Errors at time $T=1$ (using $\Delta t= 10^{-7}$) and $\alpha
    (t)$ varying in a step-wise fashion: $\alpha$ decreasing from
    $0.5$ to $0.3$ (left); $\alpha$ randomly chosen (right).}
  \label{table:T3}
\end{table}

%% ----
%% sec5
%% ----

\section{Automatic decision-making approach}
\label{sec5:automatic:decision-making:approach}

In this section, we still continue to handle the case $f=0$, whereas
the non-homogeneous case will be treated in the next section.
As we anticipated in the introductory section, we adopt two different
ML techniques in order to predict suitable values of $\alpha$, able to
guarantee stability and good performance.
When $f=0$, to train the system we take a set of functions of the
following type
\begin{align}
  \gs_{a, H}(\xs) = H \, \exp\big(-a^2\xs^2\big),
  \label{eq:train}
\end{align}
where the width $a$ of the Gaussian function is randomly chosen within
the interval $[0.2,0.6]$ and its height $H$ is also randomly chosen
within the interval $[0,1]$.
In this way we construct a total of 40 functions, that we believe to
be sufficient for our goal.
To each choice of the pair $(a,H)$, we assign the two distinct sets
\begin{align}
  \calS^{\rm PV}_{a, H}=\Big\{ \gs_{a, H}(\xiN_j), \ 1\leq j\leq N \Big\}
  \quad\textrm{and}\quad
  \calS^{\rm FC}_{a, H}=\Big\{ \widehat g_{m,a,H}, \ 0\leq m\leq N-1 \Big\},
  \label{eq:sets}
\end{align}
where the labels {\rm PV} and {\rm FC} respectively stand for {Point
  Values} and {Fourier Coefficients}.
%{\color{red} P V F C roman} 
%%
The first set clearly contains the point-values of each Gaussian
function at the Hermite nodes.
The other set is represented by the first $N$ Fourier coefficient of
$\gs_{a, H}$.
%the interpolating polynomial $\widetilde p_{N,a,H}$ of degree $N-1$ such that:
%\begin{align}
%  \widetilde\ps_{N,a,H}(\xiN_j) = g_{a,H}(\xiN_j)\exp\big( (\xiN_j)^2 \big),
%  \qquad 1\leq j\leq N.
%  \label{eq:costruzionep}
%\end{align}
%%
%Indeed, from the quadrature formula we have for $0\leq m\leq N-1$:
%\begin{align}
%  \widehat g_{m,a,H} = \frac{1}{2^m\,m!}
%  \int_{\REAL}
%  \widetilde\ps_{N,a,H} (\xs)\exp(-x^2)
%  \Hs_{m}(\xs)
%  \,dx =
%  \frac{1}{2^m\,m! \sqrt{\pi}}\sum_{j=1}^N
% g_{a,H}(\xiN_j)\exp\big( (\xiN_j)^2 \big) H_m(\xiN_j) w^N_j
%  \label{eq:coeffg}
%\end{align}
These are built as follows, for $0\leq m\leq N-1$:
%%
%% ELSEVIER
\begin{align}
  \widehat g_{m,a,H} = \frac{1}{2^m\,m!}
  \int_{\REAL}
  g_{m,a,H}(x)
  \Hs_{m}(\xs)
  \,dx \approx
  \frac{1}{2^m\,m! \sqrt{\pi}}\sum_{j=1}^N
 g_{a,H}(\xiN_j)\exp\big( (\xiN_j)^2 \big) H_m(\xiN_j) w^N_j
  \label{eq:coeffg}
\end{align}
In this way, in the special situation corresponding to $a=1$,
%one gets that $\widetilde p_{N,a,H}=H$ is a constant polynomial,
all the quantities $\widehat g_{m,a,H}$ are zero for $m\geq 1$.

The first approach is suitable to the collocation setting and the
second one to the Galerkin setting.
However, we can easily move from one setting to the other one through
the quadrature formulas of Section~\ref{sec3:quadrature}, so we can
advance the numerical solution by using the scheme
in~\eqref{eq:eulerocollocazione}, and indifferently make use of
$\calS^{\rm PV}_{a, H}$ or $\calS^{\rm FC}_{a,H}$.
To each set in~\eqref{eq:sets}, we only assign the parameter $a$, so
that Gaussian functions displaying a different height $H$ are
considered equivalent.
In this way, the algorithm is tuned on the width of the Gaussian
profiles, i.e., on what determines the choice of $\alpha$ in the
discretization of the heat equation.

%% SVM implementation
For the implementation of the SVM algorithm we use the  open source machine learning library:
\texttt{LIBSVM}~\cite{CC01a}.%\cite{Chang-Lin:2011}.
%%
%Technically,
We adopt a $\nu$-SVR ($\nu$-Support Vector Regression) approach with a radial basis function as kernel function.
The system can be either trained by using the representation suggested
in $\calS^{\rm PV}_{a, H}$ or in $\calS^{\rm FC}_{a,H}$, by
varying the parameters $a$ and $H$ in their domains.
The parameter $\nu$ allows for the control of the number of Support Vectors. The reader is
addressed for more details to~\cite{CC01a} and the references therein.
%and call the routines
%\texttt{svmtrain} and \texttt{svmpredict} as follows:

%\smallskip
%\begin{small}
%\begin{verbatim}
%   model = svmtrain(V,A,'-s 4 -t 2 -c 5.5 -e 0.0002');
%   [predictedy, accuracy, decisionvalues] = %svmpredict(V,A,model);
%\end{verbatim}
%\end{small}

\begin{figure}[h!] 
  \centering
  \includegraphics[width=1\linewidth]{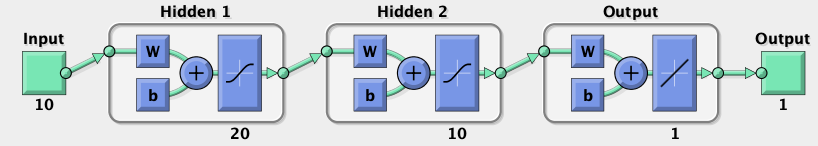}%%{./Figures/rete1.png} 
    \vskip-0.4truecm
  \caption{Matlab   neural network }
  \label{fig:reteNN} 
\end{figure}

%\smallskip
%Herein, $V$ is the matrix containing the input training values either
%in the Point-Values format or the Fourier Coefficients format, and $A$
%is the vector of the corresponding values of $a$.
%%
%This routine solves the following constrained optimization problem:
%\begin{align*}
%  \textbf{scrivere il problema di ottimizzazione.}
%\end{align*}
%Instead, for the implementation of the deep learing algorithm, se use
%the Matlab commands \texttt{fitnet} and \texttt{train} from the Matlab
%Deep Learning module.\CITE{Matlab (mettere un citazione)}
%%
%These routines allows us to build the neural network with two hidden
%layers of respectively 20 and 10 neurons each shown in Figures
%\ref{fig:reteNN}).
For the implementation of the deep learning algorithm, we use instead
the Matlab function fitting neural network: \texttt{fitnet}, which is
a feed-forward network with the tan-sigmoid transfer function in the
hidden layers and linear transfer function in the output layer.
Generally speaking, function fitting is the process of training a NN
on a set of inputs in order to produce an associated set of target
outputs.  More information can be found in the Deep Learning Toolbox
of Matlab (\texttt{www.mathworks.com}).  Based on an ``hit $\&$ trial"
approach we train with the Levenberg-Marquardt
algorithm~\cite{Kelley:1999} a NN with two hidden layers of
respectively 20 and 10 neurons (see Figure \ref{fig:reteNN}).  Note
that this network has only one output neuron, because there is only
one target value (i.e.: $a$) associated with each input vector of size
$N=10$, chosen either in $\calS^{\rm PV}_{a, H}$ or in $\calS^{\rm
  FC}_{a,H}$. %(see~\eqref{eq:sets}).

%with the following prescriptions: {\color{red} qui ho
%  buttato dentro le cose come sono state depositate da Lorella.  Vanno
%  rimesse a posto e specificate meglio. numeroDati non compare negli
%  input. Metterei anche uno schemino della rete e la regressione dopo
%  l'addestramento.}

%\smallskip
%\begin{verbatim}
%     setdemorandstream(491218382)
%     net.divideParam.trainRatio = 80/100;   
%     net.divideParam.valRatio = 10/100; 
%     net.divideParam.testRatio = 10/100;
%     neurons1=5;
%     neurons2=5;
%     numeroDati=40;  ???? (numero funzioni di training)
%     net = fitnet([neurons1 neurons2],'trainlm');
%     [net,tr] = train(net,V,A);
%\end{verbatim}

\smallskip
We then return to the experiments of the previous section by fixing
$N=10$.
We recall that the errors, when $\alpha(t)$ is taken as
in~\eqref{eq:exact:solution}, are given in the last row of Table
\ref{table:T1}.
We would like to modify $\alpha$ ten times in the time interval
$[0,1]$ according to the suggestions of the machine.
In the notation that follows the subscript {SVM} stands for
{Support Vector Machines} and subscript {DL} stands for {Deep Learning}.

 Table   \ref{table:T4}  shows the values of $\alpha$ obtained  at various times  using $\Delta t= 10^{-7}$, the two different
    representations of the training sets in~\eqref{eq:sets} and the  two     different ML algorithms.  
The errors $\displaystyle {\cal N}^1 $ and $\displaystyle {\cal N}^2 $
at time $T=1$,  $\Delta t= 10^{-7}$, for various $\alpha$, obtained with the two different
representations of the training sets in~\eqref{eq:sets} and the two
different ML techniques are shown in Table \ref{table:T5}.
These results are quite satisfactory.
The machine is actually capable to guess the appropriate value of
$\alpha$ in a given situation.
Of course, the errors improve if we decide to upgrade $\alpha$ more  frequently.

Results similar to those presented in  this section and the next one were also obtained by
trying other choices of the  network parameters, such as for example its  configuration.  
As expected,  performances may
strongly  depend upon the architecture and the  training algorithm. 
According to our experience, the most critical step is to find the right number of hidden layers 
and their respective neurons when a NN is used.

%% SIAM
% Coefficienti Fourier
\ifdefined\SiamStyle
\renewcommand{\TABROW}[6]{ \hspace{3mm} #1 \hspace{1mm} & \hspace{1mm} #2 \hspace{1mm} & \hspace{1mm} #3 \hspace{1mm} & \hspace{1mm} #4 \hspace{1mm} & \hspace{1mm} #5 \hspace{1mm} & \hspace{1mm} #6 }
\begin{table}[ht!]
  \vspace{0.3cm}
  \centering
  \begin{tabular}{c|c|cccc}
    \hline\hline
    \TABROW
        {$t$}{$\alpha =1/2 \sqrt{t+1} $}
        {$\alpha_{\rm SVM}^{\rm FC}$}
        {$\alpha_{\rm DL}^{\rm FC}$}
        {$\alpha_{\rm SVM}^{\rm PV}$}
        {$\alpha_{\rm DL}^{\rm PV}$}\\[0.5em]
        \hline\hline
    %% -------------------------------------------------------------------------------------
    \TABROW{$0$  }{ 0.5000  }{  0.5000   }{ 0.5000  }{  0.5000   }{ 0.5000   }\\
    \TABROW{$0.1$}{ 0.4767  }{  0.4736   }{ 0.4474  }{  0.4888   }{ 0.4863   }\\
    \TABROW{$0.2$}{ 0.4564  }{  0.4549   }{ 0.3639  }{  0.4691   }{ 0.4729   }\\
    \TABROW{$0.3$}{ 0.4385  }{  0.4383   }{ 0.4526  }{  0.4509   }{ 0.4613   }\\
    \TABROW{$0.4$}{ 0.4226  }{  0.4233   }{ 0.4526  }{  0.4340   }{ 0.4502   }\\
    \TABROW{$0.5$}{ 0.4082  }{  0.4097   }{ 0.4526  }{  0.4183   }{ 0.4391   }\\
    \TABROW{$0.6$}{ 0.3953  }{  0.3972   }{ 0.4526  }{  0.4039   }{ 0.4276   }\\
    \TABROW{$0.7$}{ 0.3835  }{  0.3857   }{ 0.4525  }{  0.3905   }{ 0.4159   }\\
    \TABROW{$0.8$}{ 0.3727  }{  0.3749   }{ 0.4456  }{  0.3781   }{ 0.4039   }\\
    \TABROW{$0.9$}{ 0.3627  }{  0.3649   }{ 0.3734  }{  0.3667   }{ 0.3918   }\\
    \TABROW{$1$  }{ 0.3536  }{  0.3555   }{ 0.4384  }{  0.3560   }{ 0.3800   }\\
    %% -------------------------------------------------------------------------------------
    %%\TABROW{     }{         }{           }
    \hline\hline                  
  \end{tabular}
  \vspace{0.3cm}
  \caption{Values of $\alpha$ at various times (using $\Delta t= 10^{-7}$) obtained with the two different
    representations of the training sets in~\eqref{eq:sets}, and two
    different ML algorithms.  }
  \label{table:T4}
\end{table}
\fi

%% ELSEVIER
\ifdefined\ElsevierStyle
\renewcommand{\TABROW}[6]{ \hspace{0.3cm} #1 \hspace{0.4cm} & \hspace{0.4cm} #2 \hspace{0.4cm} & \hspace{0.4cm} #3 \hspace{0.65cm} & \hspace{0.4cm} #4 \hspace{0.65cm} & \hspace{0.4cm} #5 \hspace{0.65cm} & \hspace{0.4cm} #6 }
\begin{table}[ht!]
  \vspace{0.3cm}
  \centering
  \begin{tabular}{c|c|cccc}
    \hline\hline
     \TABROW{$t$}{$\alpha =1/2 \sqrt{t+1} $}
     {$\alpha_{\rm SVM}^{\rm FC}$}  {$\alpha_{\rm DL}^{\rm FC}$}{$\alpha_{\rm SVM}^{\rm PV}$}{$\alpha_{\rm DL}^{\rm PV}$}\\[0.5em]
    \hline\hline
    \TABROW{$0$  }{ 0.5000  }{  0.5000   }{ 0.5000  }{  0.5000   }{ 0.5000   }\\[-0.5em]
    \TABROW{$0.1$}{ 0.4767  }{  0.4736   }{ 0.4474  }{  0.4888   }{ 0.4863   }\\[-0.5em]
    \TABROW{$0.2$}{ 0.4564  }{  0.4549   }{ 0.3639  }{  0.4691   }{ 0.4729   }\\[-0.5em]
    \TABROW{$0.3$}{ 0.4385  }{  0.4383   }{ 0.4526  }{  0.4509   }{ 0.4613   }\\[-0.5em]
    \TABROW{$0.4$}{ 0.4226  }{  0.4233   }{ 0.4526  }{  0.4340   }{ 0.4502   }\\[-0.5em]
    \TABROW{$0.5$}{ 0.4082  }{  0.4097   }{ 0.4526  }{  0.4183   }{ 0.4391   }\\[-0.5em]
    \TABROW{$0.6$}{ 0.3953  }{  0.3972   }{ 0.4526  }{  0.4039   }{ 0.4276   }\\[-0.5em]
    \TABROW{$0.7$}{ 0.3835  }{  0.3857   }{ 0.4525  }{  0.3905   }{ 0.4159   }\\[-0.5em]
    \TABROW{$0.8$}{ 0.3727  }{  0.3749   }{ 0.4456  }{  0.3781   }{ 0.4039   }\\[-0.5em]
    \TABROW{$0.9$}{ 0.3627  }{  0.3649   }{ 0.3734  }{  0.3667   }{ 0.3918   }\\[-0.5em]
    \TABROW{$1$  }{ 0.3536  }{  0.3555   }{ 0.4384  }{  0.3560   }{ 0.3800   }\\[-0.5em]
    \hline\hline                  
  \end{tabular}
  \vspace{0.3cm}
  \caption{Values of $\alpha$ at different times (using $\Delta t= 10^{-7}$) obtained with the two
    different representations of the training sets in~\eqref{eq:sets}
    and two different Machine Learning techniques. }
  \label{table:T4}
\end{table}
\fi

\renewcommand{\TABROW}[3]{ \hspace{0.3cm} #1 \hspace{0.4cm} & \hspace{0.4cm} #2 \hspace{0.4cm} & \hspace{0.4cm} #3 }
\begin{table}[ht!]
  \vspace{0.3cm}
  \centering
  \begin{tabular}{l|c|c}
    \hline\hline
    \TABROW{{\rm Choice of} $\alpha$}{$\displaystyle {\cal N}^1 $}{$\displaystyle {\cal N}^2$}\\
    \hline\hline
    \TABROW{ $\alpha_{\rm SVM}^{\rm FC} $}{2.6985e-08}{4.6072e-08} \\[0.5em] 
    \TABROW{ $\alpha_{\rm DL}^{\rm FC}  $}{2.9153e-05}{5.5487e-04} \\[0.5em]
    \TABROW{ $\alpha_{\rm SVM}^{\rm PV} $}{4.3810e-07}{5.9742e-07} \\[0.5em]
    \TABROW{ $\alpha_{\rm DL}^{\rm PV}  $}{6.1306e-06}{1.6049e-05} \\[0.5em]
    \hline\hline                  
  \end{tabular}
  \vspace{0.3cm}
  \caption{Errors at time $T=1$ (using $\Delta t= 10^{-7}$)  for various $\alpha$,
    obtained with the two different representations of the training
    sets in~\eqref{eq:sets} and two different ML
    techniques.}
\label{table:T5}
\end{table}

%% ----
%% sec6
%% ----

\section{The non-homogeneous case}
\label{sec6:nonhomogeneous:case}

\begin{figure}[ht!]
  \begin{center}
   \includegraphics[width=0.6\linewidth]{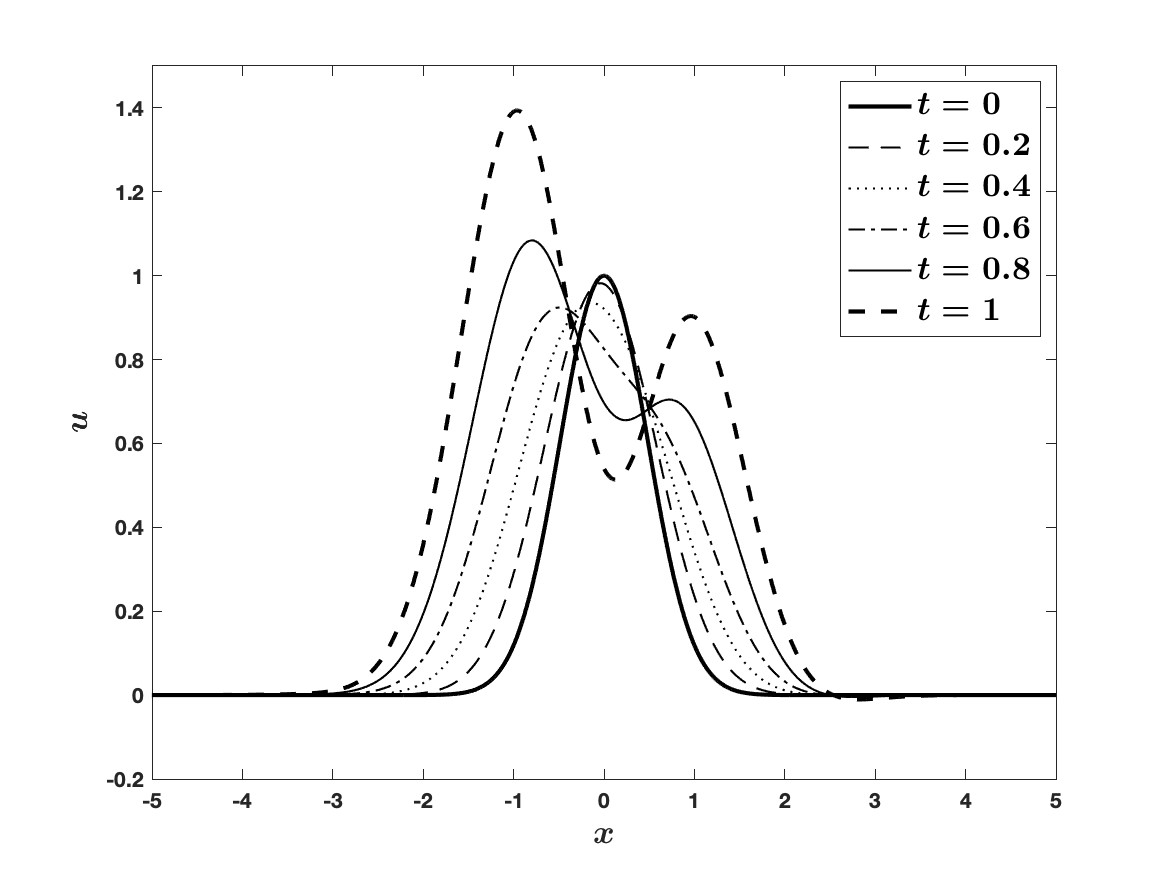}%%{./Figures/soluz.jpg} 
  \end{center}
    \vskip-0.4truecm
  \caption{The non homogeneous case: the solution $\us$ given in
    (\ref{eq:exact:solutionNO}) at different times. }
  \label{fig:solution}
\end{figure}

Based on the results obtained in the previous sections, we can now
discuss the case when $f\not =0$.
We set the right-hand side of~\eqref{eq:strong:A} in such a way that
the exact solution is:
\begin{align}
  \us(x,t) = \Big(\cos({\textstyle{\frac12}} x t) +2 (t \sin(x))^2\Big) \exp\big( -\alpha^2(t)x^2 \big)
  \quad\textrm{with}\quad
  \alpha(t) = \frac{\sqrt 2}{\sqrt{3t+1}}.
  \label{eq:exact:solutionNO}
\end{align}
With this choice $\alpha$, decays from $\sqrt{2}$ to $1/\sqrt{2}$ in
the time interval $[0,T]$ with $T=1$.
The behavior of $u$ at different times is reported in
Fig.~\ref{fig:solution}.
From a single initial bump, the solution develops by producing two
asymmetric bumps.
It is then hard to find the appropriate functional space with weight
$w_\alpha$ to formalize the theoretical problem.
From the practical viewpoint, we need the help of a ML algorithm.

As far as the approximation is concerned, from now on we set $N=16$,
that corresponds to polynomials of degree fifteen.
All the zeroes of $H_{16}$ are contained in the interval $[-5,5]$.
Moreover we choose $\Delta t= 10^{-6}$ to discretize in time.
Referring to Table \ref{tab:tab6}, we soon note that the results are
excellent when we take $\alpha$ as in~\eqref{eq:exact:solutionNO}.
Nevertheless, this expression is not known a priori.
Always according to Table \ref{tab:tab6}, the choice of $\alpha$
constantly equal to $\sqrt{2}$ (recoverable from the initial datum
$u(x,0)$) is a failure.
We also tried to guess ten values of $\alpha$ randomly chosen in the
interval [0.7,1.4].
This also turned out to be a bad choice.
The successive step is to adopt the ML strategy discussed in the
previous section.
Since in the present situation the solution $u$ is not just a Gaussian
function, we need a more representative characterization of the
training sets.
Indeed, in order to carry out our analysis, we consider the following
variant of the ML approach already experimented.

We start by defining a set of $K$ cubic splines $S_k$,
$1\leq\ks\leq\Ks$.
The support of these functions is randomly chosen in the interval
$[-c,c]$, where $c$ is a parameter that must be properly adjusted.
Such an interval is subdivided in $M+2$ equispaced break points,
including the first endpoint at $x=-c$ and the last endpoint at $x=c$,
i.e., $M$ interior break points.
The splines are required to be zero at $x=\pm c$ together with their
first derivatives, and are uniquely determined by the values attained
at the $M$ internal points.
Moreover, they are prolonged to zero outside the interval $[-c,c]$, in
order to form a set of functions with compact support in $\REAL$.
These functions are again denoted by $S_k(x)$, $1\leq\ks\leq\Ks$, for
$x\in\REAL$, and are globally $C^1(\REAL)$-regular functions.
If $\Ks\leq\Ms$, the dimension of the space generated by these splines
is less than or equal to $M$; otherwise, if $K>M$, the splines are
linearly dependent.

A given set of $K$ splines is obtained by assigning the values at the
break points in a random way.
These values must be positive and less than or equal to a given upper
bound $\calM$, which is another parameter that must properly be set in
the algorithm.
For our experiments, we choose $K=40$, $M=5$, $c=4.5$, ${\cal M}=1$.

Other sets of $K$ functions can be employed without altering the
nature and the main features of the ML algorithm, e.g., trigonometric
functions.
However, it is preferable that this family is not constituted by
polynomials, since we want to introduce elements in the training set
that are external to the space of the classical Hermite functions.

%In analogy with~\eqref{eq:costruzionep},
For each cubic spline $S_k$
%% , $1\leq k\leq K$,
and a given value of $\alpha$, we look for the interpolating
polynomial of degree less than or equal to $N-1$ such that:
\begin{align}
  \widetilde\ps_{N,k}(\xiN_j)
  = S_k(\xiN_j)\exp\big( (\alpha \xiN_j)^2 \big),
  \qquad 1\leq j\leq N.
  \label{eq:costruzionep2}
\end{align}
%function $ p_{N}(x)\exp(-\alpha^2x^2)$
Afterwards, we compute the following quantity on a finer grid:
\begin{align}
  \Gamma_{N,k}(\alpha)=
  \max_{x\in\REAL}\Big\vert S_k(x) -
   \widetilde p_{N,k}(x)\exp(-\alpha^2x^2)\Big\vert.
\end{align}
\smallskip
Assuming that the parameter $\alpha$ is varying within a given
interval $I$, we finally determine:
\begin{align}
  %  \Gamma_{N,k}(\alpha_k)=\min_{\textstyle{\alpha \in I\atop{p_N\in{\bf P}_{N}}}}G_{N,k}(\alpha ).
  \Gamma_{N,k}(\alpha_k)=\min_{\alpha \in I}\Gamma_{N,k}(\alpha ).
  \label{eq:minimum}
\end{align}
\smallskip
The Hermite function corresponding to $\alpha_k$,
\begin{align}
  g_k(x)=\widetilde p_{N,k}(x)\exp(-\alpha_k^2x^2),
  \label{eq:minimizing}
\end{align}
is called the \textit{minimizing Hermite function} associated to the
spline $S_k$. %%, $1\leq k\leq K$.
In this way, for any $S_k$, we obtains pairs like $(g_k,\alpha_k)$,
%% $1\leq k\leq K$.
%%
The above procedure extends to larger functional sets as the ones
already considered in the previous section.
In particular, if we replace $S_k$ by a Gaussian function of the type
$\gs_{a, H}$ like in \eqref{eq:up}, the \textit{minimax}-like
procedure just introduced would suggest taking $\alpha =a$ and
$\widetilde p_{N,k}$ constantly equal to $H$.

\begin{figure}[ht!]
  \centering
  \begin{tabular}{cc}
    \includegraphics[scale=0.2]{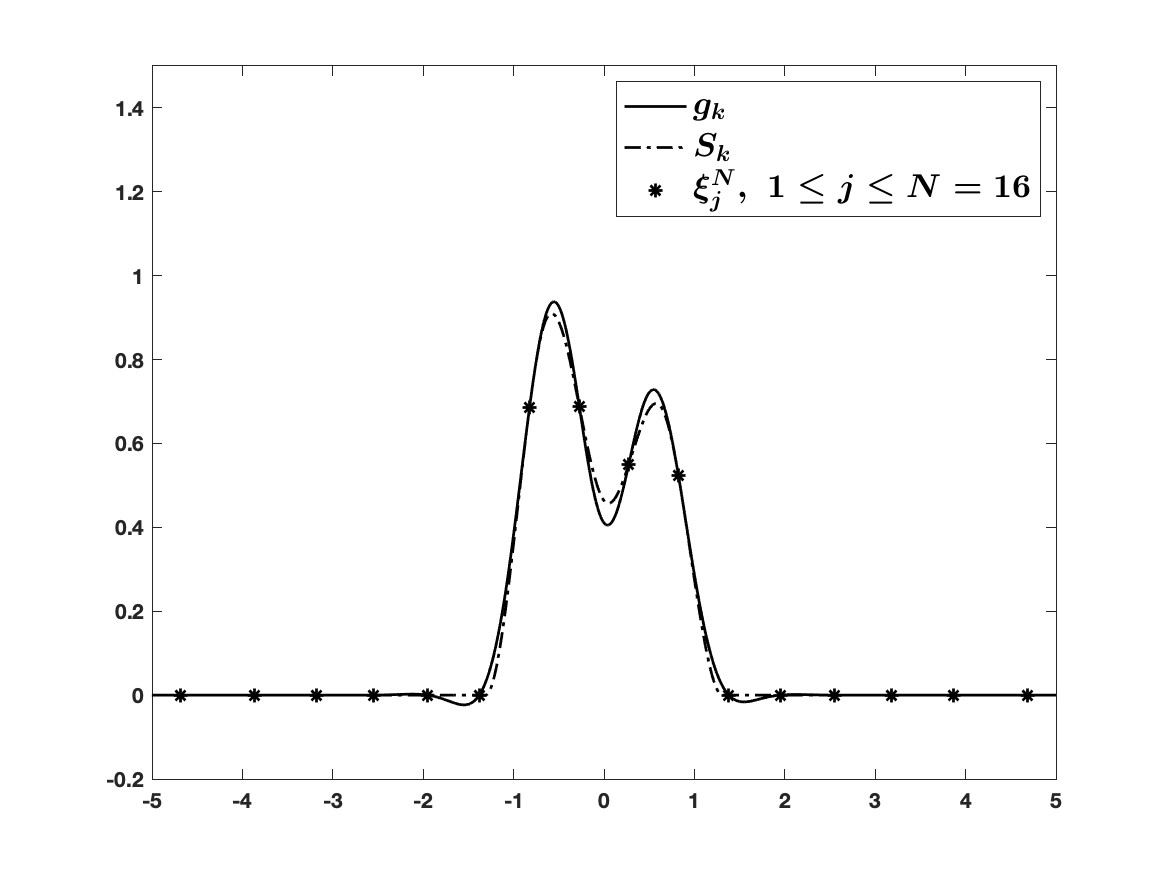} &          %% {./Figures/spline1.jpg} &
    \includegraphics[scale=0.2]{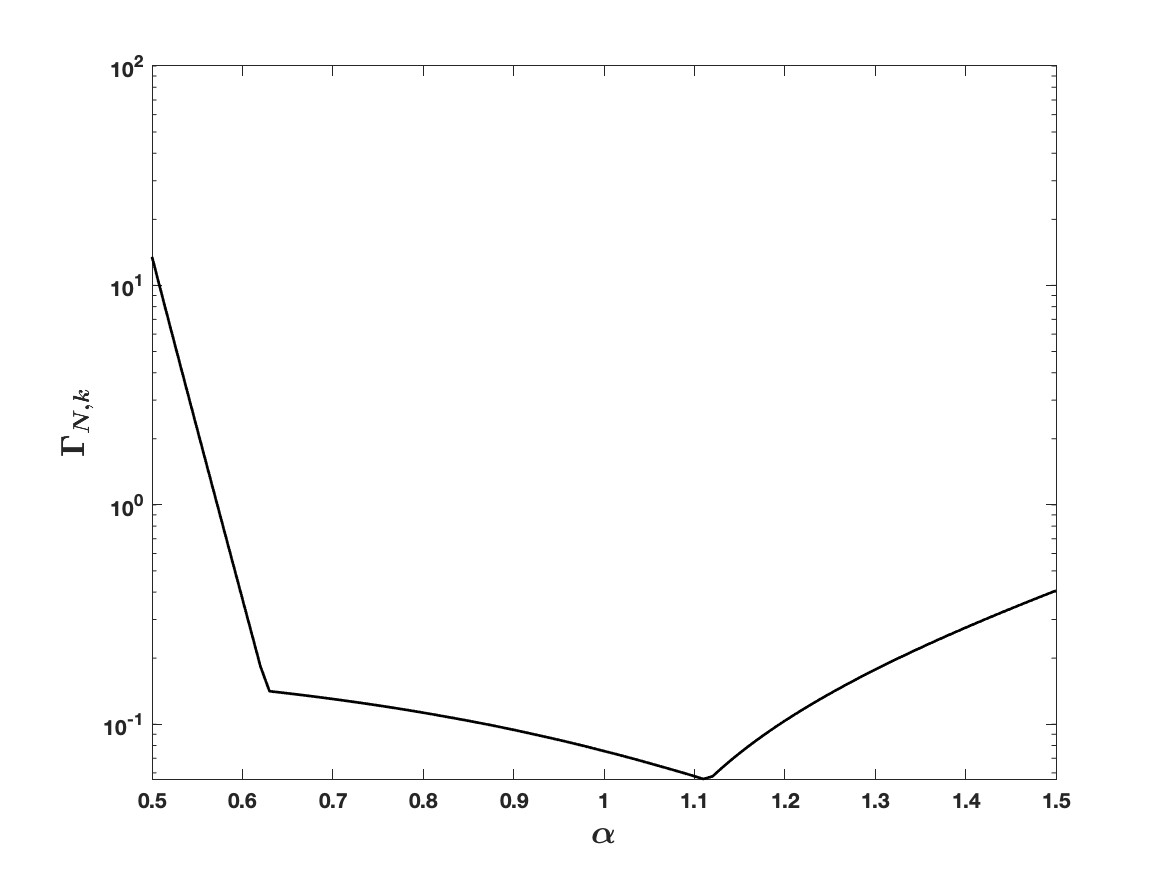} \\[0.5em]  %% {./Figures/errore1.jpg}\\[0.5em]
    %% -----------------------------------------------------------------------------
    \includegraphics[scale=0.2]{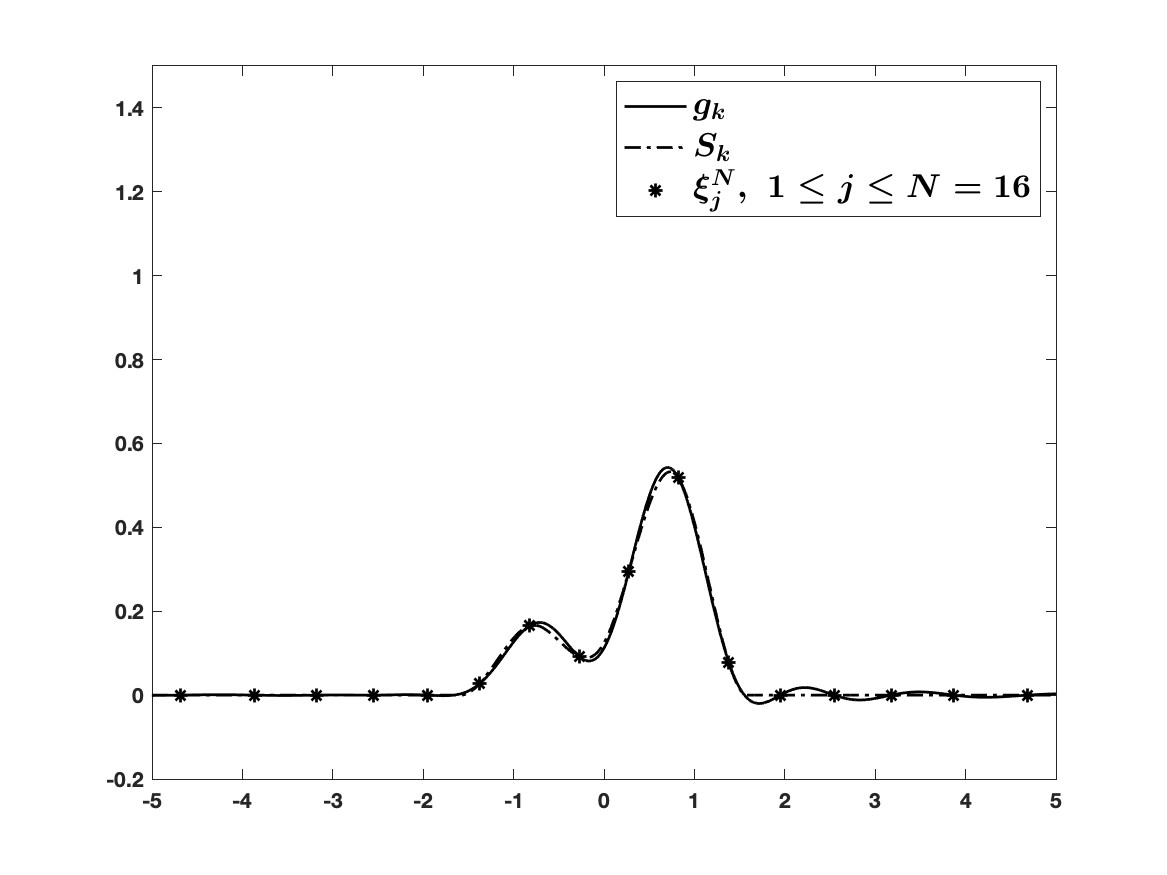} &          %% {./Figures/spline2.jpg} &
    \includegraphics[scale=0.2]{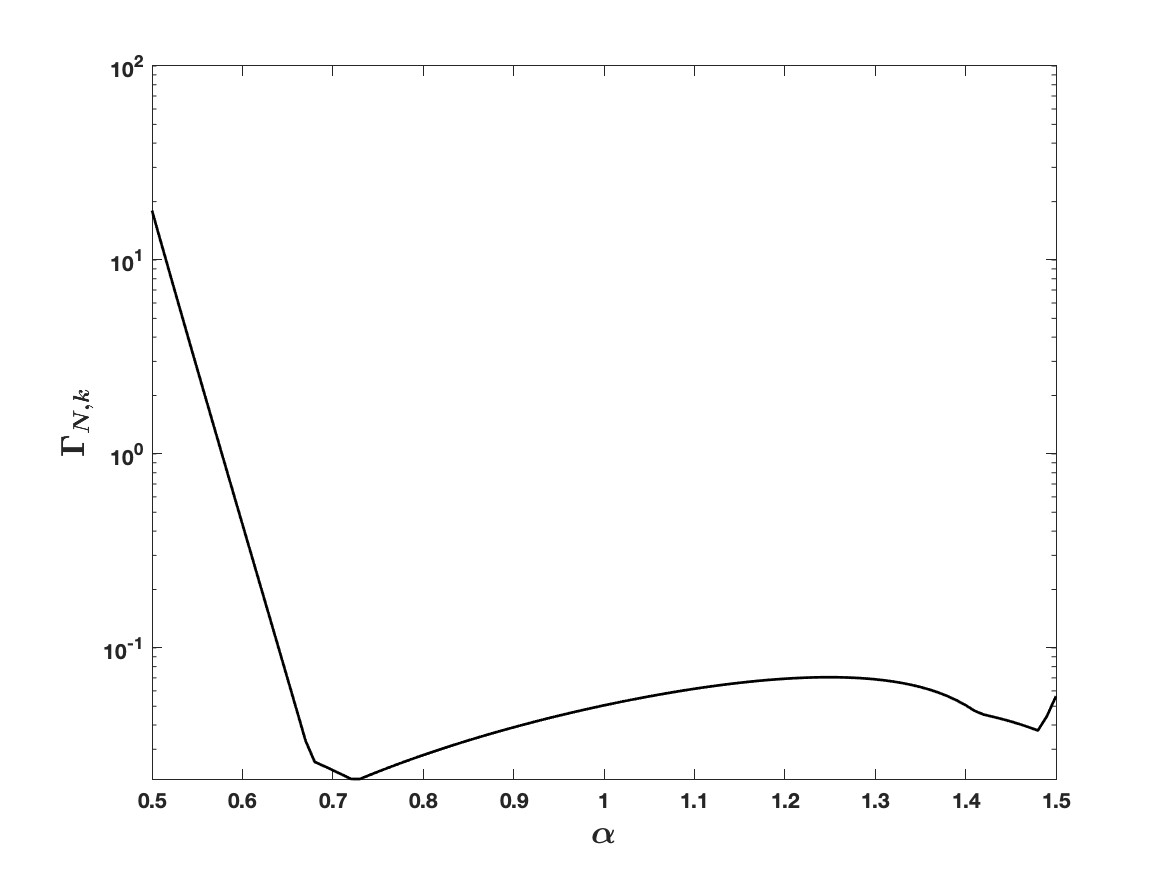}\\[0.5em]   %% {./Figures/errore2.jpg}\\[0.5em]
    %% -----------------------------------------------------------------------------
    \includegraphics[scale=0.2]{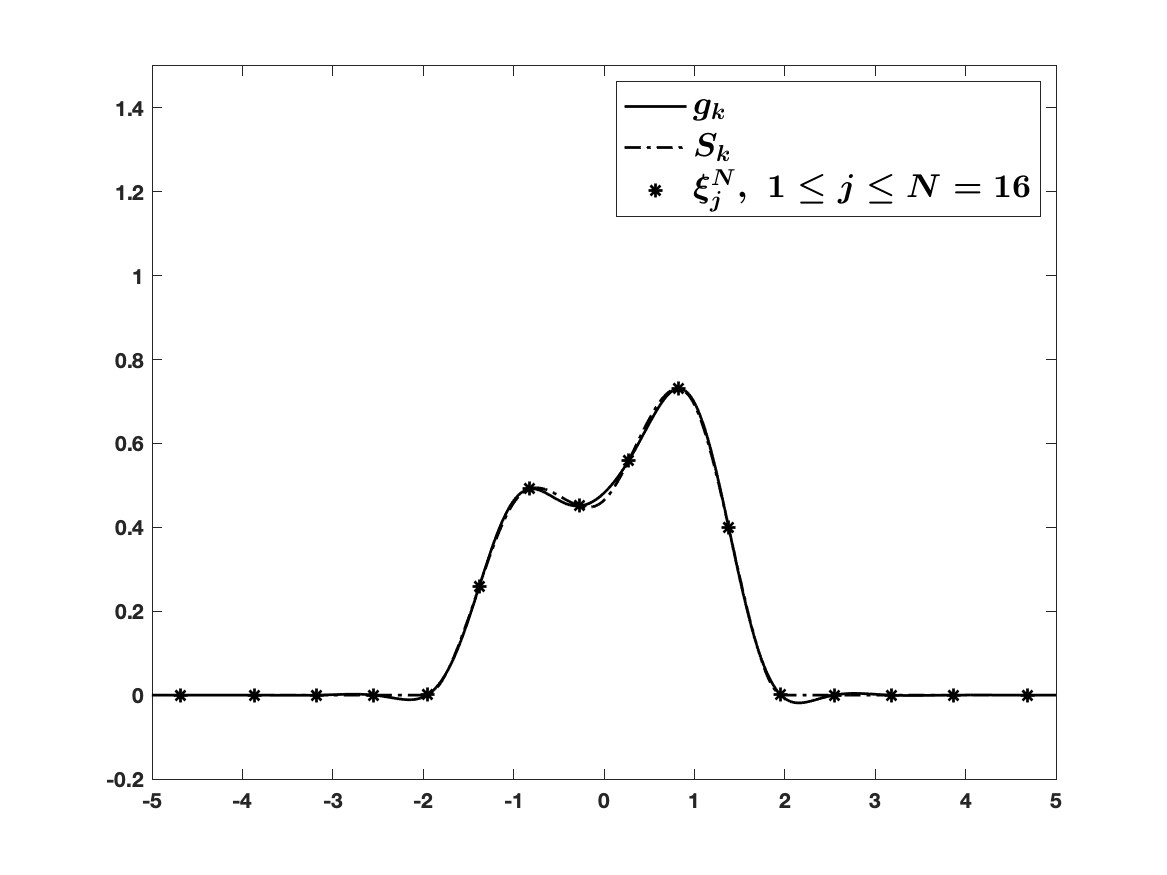} &          %% {./Figures/spline3.jpg} &
    \includegraphics[scale=0.2]{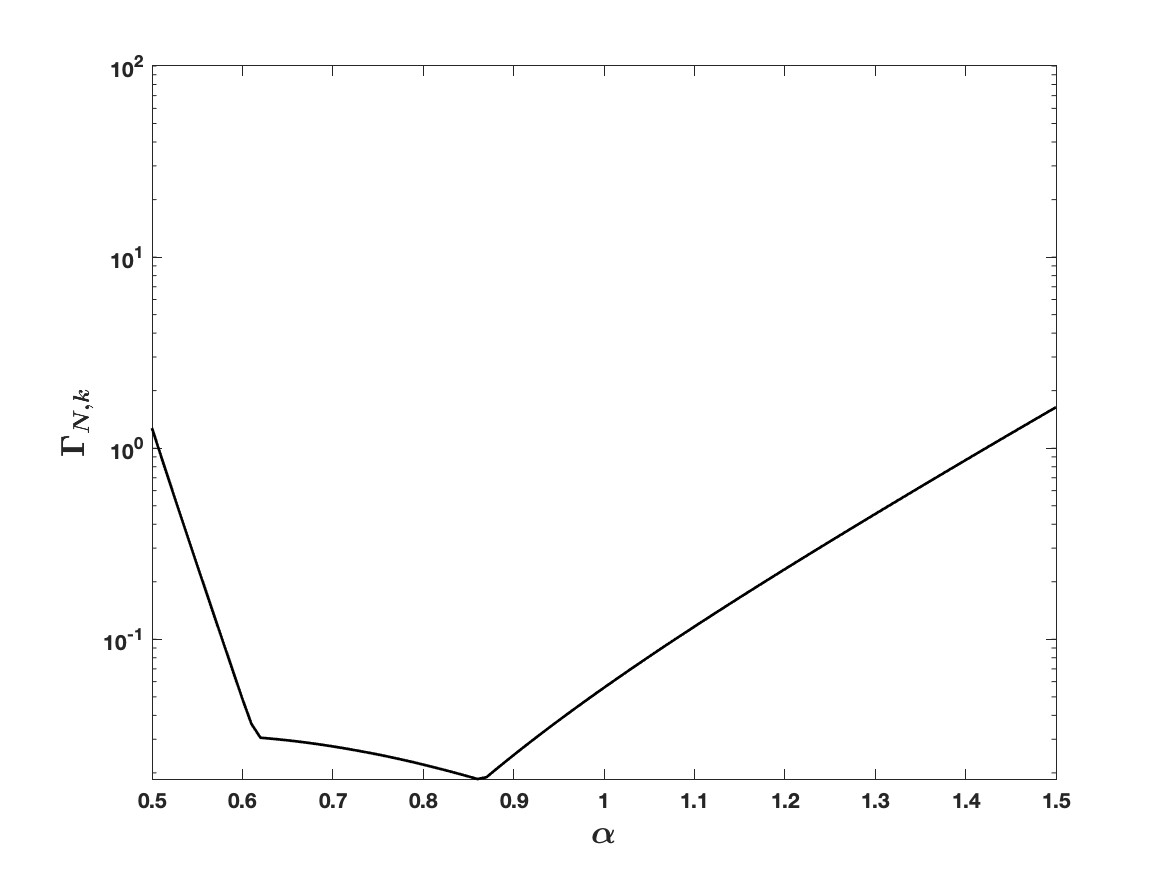} \\[0.5em]  %%{./Figures/errore3.jpg}
  \end{tabular}
  \caption{ Plots (left) of some randomly chosen splines overlapped to
    the corresponding minimizing Hermite function
    \eqref{eq:minimizing} for $x\in[-5,5]$. Plots (right) in semi-log
    scale of $\Gamma_{N,k}(\alpha)$ as a function of $\alpha$.  The
    minimizing values $\alpha_k$ are respectively 1.11 (top), 0.73
    (middle), 0.86 (bottom).  }
  \label{fig:spline}
\end{figure}

An essay of this construction is visible in Fig.~\ref{fig:spline}.
On the left, three examples of randomly generated splines are
displayed together with their corresponding minimizing Hermite
functions.
On the right, we show the semi-log plot of $\Gamma_{N,k}(\alpha)$ as a
function of $\alpha$.
It is clear that there is a minimum value $\alpha=\alpha_k$, which is
the one actually used to detect the function in \eqref{eq:minimizing}.
To obtain these graphs, we assumed that the values attained by the
splines were included in the interval $[0,{\cal M}]=[0,1]$.
Moreover, we took $I=[0.5,1.5]$.

Similarly to~\eqref{eq:sets}, we can finally define the two training
sets as follows.
The minimizing Hermite functions can be represented through the values
attained at the zeroes of the Hermite polynomial $H_{N}$.
In this way, we can build the first set of training functions:
\begin{align}
  \calS^{\rm PV}_k=\Big\{ g_k(\xiN_j)=S_k(\xiN_j), \ \  1\leq j\leq N  \Big\} , \quad\qquad  1\leq k\leq K.
  \label{eq:set2}
\end{align}
Any $\calS^{\rm PV}_k$ comes with its own $\alpha_k$, obtained
through \eqref{eq:minimum}.
The network is actually trained with the help of the values
associating $S_k$ to $\alpha_k$.
In alternative, it is possible to consider the Fourier coefficients
computed through~\eqref{eq:coeffg}.
In this case, we define the elements of a new training set as:
\begin{align}
  \calS^{\rm FC}_k=\Big\{ \widehat g_{m,k}, \ \  0\leq m\leq N-1 \Big\} , \quad\qquad  1\leq k\leq K.
  \label{eq:set3}
\end{align}

It is worthwhile to anticipate that the experiments are very sensitive
to the choice of the values of the various parameters, which are $N$,
$K$, $M$, ${\cal M}$, $c$, $I$.
In Table \ref{table:T6}, we provide the errors at $T=1$, using $\Delta
t= 10^{-6}$, for different choices of $\alpha$.
Things go smooth if we take $\alpha =\alpha (t)$ as in
\eqref{eq:exact:solutionNO}.
However this behavior is not available, with the exception of the
initial value where $\alpha (0)=\sqrt{2}$.
On the other hand, if we fix $\alpha$ to be constantly equal to
$\sqrt{2}$ the results are a disaster.
Divergence is also observed when $\alpha$ is randomly chosen in the
interval $[0.5,1.5]$ with changes at times $t^\ell =0.1\ell$, $1\leq
\ell \leq 9$.

With the automatic strategy, we again start from an initial value of
$\alpha$ equal to $\sqrt{2}$, but, then, we follow the suggestion of
the machine.
We still use the SVM algorithm and the deep learning algorithm
described in the previous section.
In particular the results of the following experiments are obtained
using a feed-forward Matlab network with two layers of five neurons
each.
As done in the previous section, the update of $\alpha$ is made 9
times in the time interval $(0,1)$.
The final results are rather good.
The parameter changes drastically at time $t=0.1$, and then, following
small variations, approximately stabilizes around the value
$\alpha\approx 0.8$ (we do not show this graph).
This is why, in the last row of Table \ref{table:T6} we report the
errors corresponding to the case where $\alpha$ has been kept fixed
for all the computation.
Without the help of the learning procedure, we would not be able to
guess that $\alpha =.8$ is a good choice.

Similar observations were made in~\cite{Fatone-Funaro-Manzini:2019b}
regarding the approximation of the 1D-1D Vlasov equation through a
periodic Fourier expansion in the space variable $x$ and Hermite
functions in the velocity variable $v$.
In that paper, for $N$ relatively small, it is numerically proven that
the impact due to a fixed choice of $\alpha$ is relevant.
The results show that the best performance is obtained for unexpected
guesses of $\alpha$, whereas the natural choice associated with the
initial datum leads to an embarrassing instability.

\renewcommand{\TABROW}[3]{ #1 & #2 & #3 }
\begin{table}[ht!]
  \vspace{0.3cm}
  \centering
  \begin{tabular}{l|c|c}
    \hline\hline
    \TABROW{{\rm Choice of} $\alpha$}{$\displaystyle {\cal N}^1 $}
           {$\displaystyle {\cal N}^2$} \\[0.5em]
    \hline\hline
    \TABROW{ $\alpha$ as in \eqref{eq:exact:solutionNO} }{ 5.8572e-04        }{ 3.2432e-04        } \\[0.3em] 
    \TABROW{ $\alpha=\sqrt{2}$                          }{ does not converge }{ does not converge } \\[0.3em] 
    \TABROW{ $\alpha$ random                            }{ does not converge }{ does not converge } \\[0.3em] 
    \TABROW{ $\alpha_{\rm SVM}^{\rm FC}$                    }{ 7.5312e-04        }{ 4.4422e-04        } \\[0.35em] 
    \TABROW{ $\alpha_{\rm DL}^{\rm FC}$                     }{ 7.5340e-04        }{ 4.5743e-04        } \\[0.35em] 
    \TABROW{ $\alpha_{\rm SVM}^{\rm PV}$                    }{ 7.5310e-04        }{ 4.4997e-04        } \\[0.35em] 
    \TABROW{ $\alpha_{\rm DL}^{\rm PV}$                     }{ 7.5330e-04        }{ 4.4862e-04        } \\[0.25em] 
    \TABROW{ $\alpha$=0.8                               }{ 7.5311e-04        }{ 4.4292e-04        } \\
    \hline\hline                  
  \end{tabular}
  \label{table:T6}
  \vspace{0.3cm}
  \caption{ Errors at time $T=1$ (using $\Delta t= 10^{-6}$) for different values of $\alpha$
    obtained with two different representations of the training sets
    (\eqref{eq:set2} and \eqref{eq:set3}) and two different ML
    techniques.}
  \label{tab:tab6}
\end{table}

One may be surprised from observing that the figures of Table
\ref{table:T6} look very similar.
We can try to give an explanation of this fact.
We realized in these numerical investigations, that the spectral
collocation method converges independently of the setting up of the
functional spaces if $\alpha$ remains within a certain range, see the
end of Section~\ref{sec2:model:problem}.
Proving convergence for a weight $w_\alpha$ that does not match
the decay of the solution $u$ is a theoretical difficulty that has
little influence on the practical outcome.
Changing $\alpha$ during the time-advancing process alters the
functional ``habitat'' of the approximated solution, but not its
capacity to produce reliable results.
Therefore, at least in the framework of the linear problem studied so
far, the guided choice of $\alpha$ is primarily useful to achieve
stability.
This is a nontrivial requisite, since the brutal choice
$\alpha=\sqrt{2}$, which agrees with the initial datum, is far from
being effective.
This phenomenon is less evident if we examine the tables relative to
the case $f=0$ (section 5), where the errors are indeed affected by
the modification of $\alpha$.
Those figures were obtained using smaller values of $N$.
In our opinion, even if the proposed values of $\alpha$ do allow
stability, there are still no sufficient nodes $\xiN_j, 1\leq j\leq
N$, in the support of the Gaussians, to guarantee an extremely
accurate approximation.
This fact agrees with the observations made in~\cite{Tang:1993} and
recalled here in the introductory section.
We also made similar observations at the end of
Section~\ref{sec3:quadrature}.
This remark stresses how critical is the set up of all the parameters
involved in such a kind of computations.

What makes the above discussion  worthwhile of attention is the fact that
to find an optimal value of the parameter $\alpha$ might not be so
crucial for improving the accuracy of the approximation, although it
could still be important for the stability.
When $N$ is sufficiently large and $\alpha$ is not exaggeratedly out of
range, the discretization method provides results that are very mildly
dependent on the choice of $\alpha$.
Of course, there are borderline situations, that can be frequent in
applications, where $N$ cannot be taken too large.
In those circumstances, an automatic technique to detect $\alpha$
comes in handy.
It has however to be said that, in general, the difficulty of finding
reliable training sets is quite a crucial issue in the organization of
ML strategy.
This may actually present a drawback in some circumstances.

We end this section with an observation.
Why did we use a ML technique instead of a plain Least Square
approach?
By assuming that the couples $\big (\{V_{k,j}\}_{1\leq j\leq N},
\alpha_k \big)$ are elements of our training set for $1\leq k\leq K$,
the standard Least Square method requires first of all the
introduction of some weights $\bar W=\{ W_j\}_{1\leq j\leq N}$.
By defining the quadratic functional
\begin{align}
  F(\bar W)=\sum_{k=1}^K\left(\sum_{j=1}^NV_{k,j}W_j \ - \ \alpha_k\right)^{\hspace{-.1cm}2},
\end{align}
we would like to find the set of weights $\bar W_{\rm min}=\{ W_{{\rm
    min},j}\}_{1\leq j\leq N}$ in such a way that
\begin{align}
  F(\bar W_{\rm min})=\min_{\bar W}F(\bar W).
  \label{eq:minpro}
\end{align}
Once this has been done, let $\{V_j^\star\}_{1\leq j\leq N}$ be a new
set of entries that correspond to some smooth function having compact
support.
We then assign a value $\alpha^\star$, according to:
\begin{align}
  \alpha^\star =\sum_{j=1}^N V_j^\star \  W_{{\rm min},j}.
  \label{eq:alpha-star:def}
\end{align}

Unfortunately, preliminary numerical results showed that such an
estimate of $\alpha^\star$ is not reliable.
The reasons for the mismatch are the following.
If $K\ll N$, there are not enough functions in the training set to
appropriately describe the new entry.
If $K>M$, there are too few break points to describe the spline set.
%%
%The ranks of the two $12\times17$-sized matrices $A^c$ and $A^p$ are
%equal to $6$ (even if we consider symmetric even splines), which is
%quite small.
%%
Thus, the splines are linearly dependent and the minimization process
to determine $\bar W$ does not give satisfactory results, since the
determinant of the Linear Regression matrix associated with the
minimization problem~\eqref{eq:minpro} turns out to be zero.
If $M\approx N$, there are too many break points.
For this reason, the splines may present oscillations which are not
properly caught by the corresponding minimizing Hermite function (the
situation is even worse if $M>N$).
In all cases, the learning procedure may be considered unsatisfactory.
Instead, a nonlinear ML approach allows us to introduce a larger
number of weights in the system.
As a consequence, we can play with more functions in the training set,
and this is true also when $K>M$, i.e., when the splines start to be
linearly dependent.
The just mentioned situation is rather different from those emerging, for instance,
from problems in Image Recognition.
%%\CITE{at least one paper!}{\color{red}Fuck! Non lo so}
%%
There, the number $K$ of pictures in the training set is much less
that $M$, which is the number of pixels needed to represent a single
image.

%% -----
%% sec 7
%% -----

\section{Conclusions and future work}
\label{sec7:conclusions}

In this work, we proposed an automatic decision-making system for the
determination of the scaling factor $\alpha$ of the weight function in
the Hermite spectral method for a PDEs on unbounded domains.
An appropriate value of such parameter is crucial when using this kind
of approximations since a bad choice may result in instabilities or
impractical costs of implementation.
We employed ML techniques based on either deep NN or SVM, in order to
predict $\alpha$ in automatic way, in the case of a one-dimensional
time-dependent equation.
After the training of the machine, the algorithm advances in time by
updating periodically $\alpha$ according to a procedure that extracts
the value of the parameter in agreement to the behavior manifested in
phase of evolution by the solution itself.
The numerical investigations carried out show that the algorithm is
able to determine $\alpha$ and maintain the stability, so improving
when possible the accuracy of the outcome.
As a proof of concept and for the exposition sake, we applied this
approach to the 1D heat equation, but we are confident that
applications to
%% \st{\textbf{more serious}} 
non-linear multi-dimensional problems can turn out to be successful,
in particular in the domain of plasma physics.

In the case of the 1D-1V Vlasov equation, the best known algorithm to
select the shifting and scaling factors relies on a physically-based
criterion, by following the average velocity and temperature of each
plasma species, i.e., the first and second moments of the plasma phase
space density~\cite{Delzanno-Manzini-Pagliantini-Markidis:2019}.
This time-dependent adaptive strategy has been proved to preserve the
conservation laws of total mass, momentum and energy of the
non-adaptive approach.
All these considerations dictate the direction of our future work,
i.e., the design of ML strategies for the automatic determination of
both $\alpha$ and $\beta$ in the spectral approximations of the Vlasov
equation. This fulfillment requires additional efforts, due to the
fact that $\alpha =\alpha(t, \vec x)$ and $\beta=\beta(t, \vec x)$ now
also depend on the location of the particles.
A straightforward implication is that the design of the training set
becomes much more complicated.

%% Acknowledgments
%% \input{acknowledgments.tex}
\section*{Acknowledgments}
The Authors are grateful to Dr. G.~L. Delzanno (LANL) 
and Prof. C. Pagliantini for many fruitful discussions
and suggestions.
The Authors are affiliated to GNCS-INdAM (Italy).
The third author was supported by the LDRD program of Los
Alamos National Laboratory under project number 20170207ER.
Los Alamos National Laboratory is operated by Triad National Security,
LLC, for the National Nuclear Security Administration of
U.S. Department of Energy (Contract No. 89233218CNA000001).

%% BIBLIOGRAPHY
%\clearpage
%\bibliographystyle{plain}
%\bibliography{hermite}

%% APPENDICES
\clearpage
%%\input{new_appendix.tex}
% Hey Emacs, this is -*-latex-*-

%\appendix

\section*{Appendix}
%\subsection{Hermite polynomials}
We begin by collecting some basic relations concerning Hermite
polynomials.
We denote with a prime the derivative of a given  Hermite
polynomial with respect to the argument $\zeta$.
Thus, the first and the second derivatives of the $\ell$-th Hermite polynomial, for $\ell \geq 2$, are given by:
%\begin{align}
%  \Hsp_{\ell}(\zeta) =
%  \begin{cases}
%    0                      & \ell=0,\\
%    2\ell\Hs_{\ell-1}(\zeta) & \ell>0.
%  \end{cases}\qquad\qquad
%  \begin{cases}
%    0                                           & \ell<2,\\
%    4\displaystyle{\frac{\ell!}{(\ell-2)!}}\Hs_{\ell-2}(\zeta) & \ell\geq2.
%  \end{cases}
%  \label{eq:app:Hermite:first:derivative}
%\end{align}
%Alternatively, we can write:
\begin{align}
 \Hsp_{\ell}(\zeta)=2\ell\Hs_{\ell-1}(\zeta),\qquad\qquad
  \Hspp_{\ell}(\zeta) = 4\ell(\ell-1)\Hs_{\ell-2}(\zeta).
   \label{eq:app:Hermite:derivatives}
\end{align}
In addition, we have: $H'_0(\zeta)=0$, $H'_1(\zeta)=2$, $H''_0(\zeta)=0$, $H''_1(\zeta)=0$,
so that~\eqref{eq:app:Hermite:derivatives} formally holds for any $\ell\geq 0$.
%for $\ell\geq2$ and note that the second derivative of
%$\Hs_{0}(\zeta)=1$ and $\Hs_{1}(\zeta)=2\zeta$ is obviously zero.
%%
Useful recursive relations are:
\begin{align}
  2\zeta\Hsp_{\ell}(\zeta)
  = \Hspp_{\ell}(\zeta) + 2\ell\Hs_{\ell}(\zeta)
  =
  \begin{cases}
    2\ell\Hs_{\ell}(\zeta) ,                                 & \ell<2,\\
    4\ell(\ell-1)\Hs_{\ell-2}(\zeta) + 2\ell\Hs_{\ell}(\zeta) ,& \ell\geq2.
  \end{cases}
  \label{eq:app:Hermite:useful:relation}
\end{align}
%%
%Furthermore, reversing~\eqref{eq:app:Hermite:first:derivative} for
%$\ell>0$, we find that
%$\Hs_{\ell-1}(\zeta)=\Hsp_{\ell}(\zeta)\slash{(2\ell)}$.
%%
%Translating $\ell\to\ell+1$, we find that 
%\begin{align}
%  \Hs_{\ell}(\zeta) = \frac{1}{2(\ell+1)}\Hsp_{\ell+1}(\zeta)
%  \qquad\ell>0.
%\end{align}
%Then, we multiply both sizes by $\zeta$; we
%use~\eqref{eq:app:Hermite:useful:relation} (again translated to
%$\ell+1$) and find that
\begin{align}
  \zeta\Hs_{\ell}(\zeta)
  &= \frac{1}{4(\ell+1)}\,2\zeta\Hsp_{\ell+1}(\zeta)
  = \frac{1}{4(\ell+1)}\Big[4(\ell+1)\ell\Hs_{\ell-1}(\zeta) + 2(\ell+1)\Hs_{\ell+1}(\zeta)\Big]
  \nonumber\\[0.5em]
  &= \ell\Hs_{\ell-1}(\zeta) + \frac{1}{2}\Hs_{\ell+1}(\zeta) ,
  \qquad\ell\geq1.   \label{eq:app:zeta:Hl}
\end{align}
For completeness, we note that for $\ell=0$ it holds
$\zeta\Hs_{0}(\zeta)=\Hs_{1}(\zeta)/2$.

We now go through the computations relative to the three formulations
\eqref{eq:weak:form:IIb}, \eqref{eq:weak:form:Ia}, \eqref{eq:cinesi}.
We start with the second one.

We recall that $\alpha$ depends on $t$ and we denote by $\alpha'$ its derivative. Concerning Hermite polynomials, the prime will continue to
denote the derivative with respect to $\zeta$. We substitute the definitions~\eqref{eq:test:function}, ~\eqref{eq:trial:function}   of, respectively,  $\phi_{\ms}$,   $\usN$, and we split the integral in three
parts that will be computed separately:
\begin{align}
  &\int_{\REAL}\pt\usN\,\phi_{\ms}\dx
  =
  \int_{\REAL}
  \frac{\pt\wsa(x,t)}{\sqrt{\pi}}\left[\sum_{\ell=0}^{N}\ush_{\ell}(t)\Hs_{\ell}\big(\alpha\xs\big)\right]
  \frac{\alpha}{2^m\,m!}\Hs_{m}\big(\alpha\xs\big)\dx
  \nonumber\\[0.5em]
  &\qquad\phantom{=}
  +
  \int_{\REAL}
  \frac{\wsa(x,t)}{\sqrt{\pi}}\sum_{\ell=0}^{N}\bigg[
    \Big(\pt\ush_{\ell}(t)\Big)\Hs_{\ell}\big(\alpha\xs\big) +
    \ush_{\ell}(t)\,\alpha'\xs\Hsp_{\ell}\big(\alpha\xs\big)
    \bigg]
  \frac{\alpha}{2^m\,m!}\Hs_{m}\big(\alpha\xs\big)\dx  
  \nonumber\\[0.5em]
  &\qquad
  =
  \frac{\alpha}{2^m\,\ms!\sqrt{\pi}}(-2\alpha\asp)
  \sum_{\ell=0}^{N}\ush_{\ell}(\ts)\int_{\REAL}\xs^2\wsa(\xs,\ts)\Hs_{\ell}(\alpha\xs)\Hs_{m}(\alpha\xs)\dx
  \nonumber\\[0.5em]
  &\qquad\phantom{=}
  +
  \frac{1}{2^{m}\,\ms!\,\sqrt{\pi}}
  \sum_{\ell=0}^{N}\left[
    \pt\ush_{\ell}(t)\int_{\REAL}\wsa(x,t)\Hs_{\ell}\big(\alpha\xs\big)\Hs_{m}\big(\alpha\xs\big)\,\alpha\dx
    \right.
    \nonumber\\[0.5em]
    &\qquad
   % =\frac{1}{2^{m}\,\ms!\,\sqrt{\pi}}\sum_{\ell=0}^{N}\left[\right.
  \left.
  +\ush_{\ell}\asp\int_{\REAL}\wsa(x,t)\,\xs\Hsp_{\ell}\big(\alpha\xs\big)\,\Hs_{m}\big(\alpha\xs\big)\,\alpha\dx
  \right]
  = \TERM{I}{} + \TERM{II}{} + \TERM{III}{}. 
\end{align}
To compute \TERM{I}{}, we substitute $\zeta=\alpha\xs$ and use
the orthogonality properties of the Hermite polynomials to find
that:
\begin{align}
  \TERM{I}{}
  &=
  \frac{1}{2^m\,\ms!\sqrt{\pi}}\frac{-2\asp}{\alpha}
  \sum_{\ell=0}^{N}\ush_{\ell}(\ts)\int_{\REAL}(\alpha\xs)^2\exp{\big(-(\alpha\xs)^2\big)}\Hs_{\ell}(\alpha\xs)\Hs_{m}(\alpha\xs)\alpha\dx
  \nonumber\\[0.5em]
  &=
  -\frac{\asp}{\alpha}
  \Big[(2\ms+1)\,\ush_{\ms}(\ts)+2(\ms+2)(\ms+1)\,\ush_{\ms+2}(\ts)+\frac{1}{2}\ush_{\ms-2}(\ts)\Big].\nonumber
   \end{align}
  %% \medskip
  %% \noindent
  %% %%
 To compute the term \TERM{II}{}, we substitute
    $\zeta=\alpha\xs$ and use the orthogonality
    properties:
  \begin{align}
  \TERM{II}{}
  &= \frac{1}{2^{m}\,\ms!\,\sqrt{\pi}}\sum_{\ell=0}^{N}\pt\ush_{\ell}(t)
  \int_{\REAL}\exp(-\zeta^2)\Hs_{\ell}(\zeta)\Hs_{m}(\zeta)\,d\zeta
  = \pt\ush_{\ms}(\ts) .\nonumber
  %%\end{align}
  \intertext{To compute term \TERM{III}{} we use~\eqref{eq:app:Hermite:derivatives}
  and the orthogonality of  Hermite polynomials:}
  %%\begin{align}
  \TERM{III}{}
  &= \frac{1}{2^{m}\,\ms!\,\sqrt{\pi}}\sum_{\ell=0}^{N}\ush_{\ell}(t)\frac{\asp}{\alpha}
  \int_{\REAL}\exp(-\zeta^2)\,\zeta\Hsp_{\ell}(\zeta)\Hs_{m}(\zeta)\,d\zeta\nonumber\\[0.5em]
  &= \frac{1}{2^{m}\,\ms!}\sum_{\ell=0}^{N}\ush_{\ell}(t)\frac{\asp}{\alpha}\Big[
    2\ell\ms\,2^{\ms-1}(\ms-1)!\,\delta_{\ell-1,\ms-1}
    + \ell   \,2^{\ms+1}(\ms+1)!\,\delta_{\ell-1,\ms+1}
    \Big]
  \nonumber\\[0.5em]
  &=
  \frac{\asp}{\alpha}\Big[
    \ms\ush_{\ms}(\ts) + 
    2(\ms+2)(\ms+1)\ush_{\ms+2}(\ts)
    \Big]. \nonumber
\end{align}
By collecting the results for $\TERM{I}{}$, $\TERM{II}{}$, and
$\TERM{III}{}$, we find out that:
\begin{align}
  \TERM{I}{} + \TERM{II}{} + \TERM{III}{}
  =
  \pt\ush_{\ms}(\ts) +
  \frac{\asp}{\alpha}\Big[
    -(\ms+1)\ush_{\ms}(\ts) - \frac{1}{2}\ush_{\ms-2}(\ts)
    \Big].
  \label{eq:app:formI:intg:time:derivative}
\end{align}
In a similar fashion, we split the integral of the second derivative of $\usN$
against the test function into three parts:
\begin{align*}
  &\int_{\REAL}\pxx\usN\phi_{\ms}\dx
  = \int_{\REAL}\pxx\wsa(\xs,\ts)
  \bigg[\frac{1}{\sqrt{\pi}}\sum_{\ell=0}^{N}\ush_{\ell}(\ts)\Hs_{\ell}(\alpha\xs)\bigg]\,
  \bigg[\frac{\alpha}{2^m\,m!}\Hs_{m}\big(\alpha\xs\big)\bigg]\dx
  \nonumber\\[0.5em]
  &\qquad\qquad
  + \frac{1}{\sqrt{\pi}}
  \int_{\REAL}
  \bigg[2\px\wsa(\xs,\ts)\sum_{\ell=0}^{N}\ush_{\ell}(\ts)\px\Hs_{\ell}(\alpha\xs)\bigg]
  \bigg[\frac{\alpha}{2^m\,m!}\Hs_{m}\big(\alpha\xs\big)\bigg]
  \dx
  \nonumber\\[0.5em]
  &\qquad\qquad
  +
   \frac{1}{\sqrt{\pi}}\frac{\alpha}{2^m\,m!}
  \int_{\REAL}
  \bigg[\wsa(\xs,\ts)\sum_{\ell=0}^{N}\ush_{\ell}(\ts)\pxx\Hs_{\ell}(\alpha\xs)\bigg]\,
  \bigg[\Hs_{m}\big(\alpha\xs\big)\bigg]\dx
  \nonumber\\[0.5em]
  &\qquad\qquad
  =
  \TERM{A}{} + \TERM{B}{} + \TERM{C}{}.
\end{align*}
Using the property of the Hermite polynomials, we evaluate these terms
as follows:
\begin{align*}
  \TERM{A}{} &= 2\alpha^2\bigg[
    2\ms\ush_{\ms}(\ts)
    + 2(\ms+2)(\ms+1)\ush_{\ms+2}(\ts)
    +\frac{1}{2}\ush_{\ms-2}(\ts)
    \bigg]
  \\[0.5em]
  %% -----------------------------
  \TERM{B}{} &= -4\alpha^2\Big[\ms\ush_{\ms}(\ts) + 2(\ms+2)(\ms+1)\ush_{\ms+2}(\ts)\Big]
  \\[0.5em]
  %% -----------------------------
  \TERM{C}{} &= 4(\ms+2)(\ms+1)\alpha^2\ush_{\ms+2}(\ts).
\end{align*}
Putting all together we arrive at:     
\begin{align}
  \TERM{A}{}+\TERM{B}{}+\TERM{C}{}
  &=
  2\alpha^2\Big[ 2\ms\ush_{\ms}(\ts) + 2(\ms+2)(\ms+1)\ush_{\ms+2}(\ts) +\frac{1}{2}\ush_{\ms-2}(\ts) \Big] \nonumber\\[0.5em]
  & -4\alpha^2\Big[\ms\ush_{\ms}(\ts) + 2(\ms+2)(\ms+1)\ush_{\ms+2}(\ts)\Big]               \nonumber\\[0.5em]
  & +4(\ms+2)(\ms+1)\alpha^2\ush_{\ms+2}(\ts)                                             = \alpha^2\ush_{\ms-2}(\ts).
  \label{eq:app:formI:intg:abc}
\end{align}
By equating \eqref{eq:app:formI:intg:time:derivative} and
\eqref{eq:app:formI:intg:abc} we finally obtain the scheme \eqref{eq:form0:final:15}.
%The variational formulation \TERM{W1}{} reads as:
%%
%\begin{align}
%  \pt\ush_{\ms}(\ts) +
%  \frac{\asp(\ts)}{\alpha(\ts)}\Big[
%    -(\ms+1)\ush_{\ms}(\ts) - \frac{1}{2}\ush_{\ms-2}(\ts)
%    \Big]
%  -\nu\alpha(\ts)^2\ush_{\ms-2}(\ts) = \fs_{\ms}(\ts),
%  \label{eq:app:form0:final:10}
%\end{align}
%and after a straightforward calculation, we find that:
%\begin{align}
%  \pt\ush_{\ms}(\ts) -
%  \bigg(\frac{\asp(\ts)}{2\alpha(\ts)} + \nu\alpha(\ts)^2\bigg)\ush_{\ms-2}(\ts)
%  -(\ms+1)\frac{\asp(\ts)}{\alpha(\ts)}\ush_{\ms}(\ts)
%  = \fs_{\ms}(\ts) 
%  \label{eq:app:form0:final:15}
%\end{align}

%\subsection{Variational formulation~\TERM{W2}{}}
We then examine the scheme originating from \eqref{eq:weak:form:IIb}.
%We substitute the definition of $\usN$ from~\eqref{eq:trial:function}
%and $\phi_{\ms}$ from~\eqref{eq:test:function} apply chain rule
%derivation:
%%
We must compute:
\begin{align*}
  & \int_{\REAL}\frac{\partial\usN}{\partial\xs}\,\frac{\partial\phi_{\ms}}{\partial\xs}\dx
  =
  \frac{1}{\sqrt{\pi}}\,\frac{\alpha}{2^{\ms}\,\ms!}
  \sum_{\ell=0}^{N}\ush_{\ell}(\ts)
  \int_{\REAL}\px\Big(\wsa(\xs,\ts)\Hs_{\ell}(\alpha\xs)\Big)\,\px\Hs_{\ms}(\alpha\xs)\,\dx
\end{align*}

\begin{align}
  &\qquad=
  \frac{1}{\sqrt{\pi}}\,\frac{\alpha}{2^{\ms}\,\ms!}\sum_{\ell=0}^{N}\ush_{\ell}(\ts)\bigg[
    \int_{\REAL}\big(\px\wsa(\xs,\ts)\big)\Hs_{\ell}(\alpha\xs)\,\px\Hs_{\ms}(\alpha\xs)\,\dx
    \nonumber\\[0.5em]
    &\qquad\phantom{=\frac{1}{\sqrt{\pi}}\,\frac{\alpha}{2^{\ms}\,\ms!}\sum_{\ell=0}^{N}\ush_{\ell}(\ts)\bigg[}
    +\int_{\REAL}\wsa(\xs,\ts)\px\Hs_{\ell}\big(\alpha\xs)\,\px\Hs_{\ms}\big(\alpha\xs\big)\,\dx
    \bigg]
  \nonumber
   \end{align}
  \begin{align}
  &\qquad=
  \frac{1}{\sqrt{\pi}}\,\frac{\alpha}{2^{\ms}\,\ms!}\sum_{\ell=0}^{N}\ush_{\ell}(\ts)\bigg[
    -2\alpha\int_{\REAL}\big(\alpha\xs\big)\wsa(\xs,\ts)\Hs_{\ell}(\alpha\xs)\,2\ms\Hs_{\ms-1}(\alpha\xs)\,\alpha\dx
    \nonumber\\[0.5em]
    &\qquad\phantom{=\frac{1}{\sqrt{\pi}}\,\frac{\alpha}{2^{\ms}\,\ms!}\sum_{\ell=0}^{N}\ush_{\ell}(\ts)\bigg[}
    +\alpha\int_{\REAL}\wsa(\xs,\ts)\Hsp_{\ell}(\alpha\xs)\,\Hsp_{\ms}(\alpha\xs)\,\alpha\dx
    \bigg]
  %  &\hspace{-2.5cm}\mbox{\big[set $\zeta=\alpha(\ts)\xs$\big]} 
  \nonumber\\[0.5em]
  &\qquad=
  \frac{1}{\sqrt{\pi}}\,\frac{\alpha^2}{2^{\ms}\,\ms!}\sum_{\ell=0}^{N}\ush_{\ell}(\ts)\bigg[
    -4\ms\int_{\REAL}\zeta\Hs_{\ell}(\zeta)\Hs_{\ms-1}(\zeta)e^{-\zeta^2}\,d\zeta
    \nonumber\\[0.5em]
  &\qquad\phantom{\hspace{4.cm}} +\int_{\REAL}\exp{(-\zeta^2)}\,\Hsp_{\ell}(\zeta)\,\Hsp_{\ms}(\zeta)e^{-\zeta^2}\,d\zeta\bigg]
  \nonumber\\[0.5em]
  &\qquad
  = \alpha^2\Big( -2\ms\,\ush_{\ms}(\ts)-\,\ush_{\ms-2}(\ts) \Big)
  + \alpha^2\Big( 2\ms\,\ush_{\ms}(\ts) \Big)
  = - \alpha^2 \ush_{\ms-2}(\ts).
  \label{eq:app:formI:intg:space:derivative:X}
\end{align}
By changing the sign of the last term in 
\eqref{eq:app:formI:intg:space:derivative:X} we exactly get the
same result in
\eqref{eq:app:formI:intg:abc} which brings again to the scheme \eqref{eq:form0:final:15}.

We finally examine the scheme originating from \eqref{eq:cinesi}.
The first integral is clearly equal to $\ush_{\ms}(\ts)$.
Successively, we evaluate:
\begin{align}
  %% ===============
  %% second integral
  %% ===============
  \int_{\REAL}\usN\,\frac{\partial^2\phi_{\ms}}{\partial\xs^2}\,\dx
  &=\alpha^2\int_{\REAL}\usN\,\phi_{\ms-2}\,\dx
 % &\hspace{-3.cm}\mbox{\big[use~\eqref{eq:trial:function} and~\eqref{eq:test:function}\big]}
  \nonumber\\[0.5em]
  & = \alpha^2\frac{1}{\sqrt{\pi}}\,\frac{\alpha}{2^{\ms-2}(\ms-2)!}\sum_{\ell=0}^{N}\ush_{\ell}(\ts)
  \int_{\REAL}\wsa(\xs,\ts)\Hs_{\ell}(\alpha\xs)\,\Hs_{\ms-2}(\alpha\xs)\,\dx
 % &\hspace{-3.cm}\mbox{\big[set $\zeta=\alpha(\ts)\xs$\big]}
  \nonumber\\[0.5em]
  &= \frac{\alpha^2}{\sqrt{\pi}}\frac{1}{2^{\ms-2}(\ms-2)!}\sum_{\ell=0}^{N}\ush_{\ell}(\ts)
  \int_{\REAL}\Hs_{\ell}(\zeta)\Hs_{\ms-2}(\zeta)e^{-\zeta^2}\,d\zeta
  % &\hspace{-3.cm}\mbox{\big[set \eqref{eq:ortho:intg:00}\big]}
  \nonumber\\[0.5em]
  &= \frac{\alpha^2}{2^{\ms-2}(\ms-2)!\sqrt{\pi}}\sum_{\ell=0}^{N}\ush_{\ell}(\ts)\,2^{\ms-2}(\ms-2)!\sqrt{\pi}\delta_{\ell,\ms-2}
  = \alpha^2\ush_{\ms-2}(\ts).\nonumber
\end{align}
Regarding the last integral, we  split it into two parts:
\begin{align}
  %% ==============
  %% third integral
  %% ==============
  &\int_{\REAL}\usN\,\frac{\partial\phi_{\ms}}{\partial\ts}\,\dx
  =
  \frac{1}{\sqrt{\pi}}\,\frac{1}{2^{\ms}\,\ms!}\sum_{\ell=0}^{N}\ush_{\ell}(\ts)
  \int_{\REAL}\wsa(\xs,\ts)\Hs_{\ell}(\alpha\xs)\,\pt\Big(\alpha\Hs_{\ms}(\alpha\xs)\Big)\dx
  % &\mbox{\big[use derivative chain rule\big]}
  \nonumber\\[0.5em]
  &\qquad
  = %%\hspace{-1cm}=
  \frac{1}{\sqrt{\pi}}\,\frac{1}{2^{\ms}\,\ms!}\sum_{\ell=0}^{N}\ush_{\ell}(\ts)
  \int_{\REAL}\wsa(\xs,\ts)\Hs_{\ell}(\alpha\xs)\,
  \Big(\asp\Hs_{\ms}(\alpha\xs)+\alpha\asp\xs\Hsp_{\ms}(\alpha\xs)\Big)\dx
  \nonumber\\[0.5em]
  &\qquad
  = %%\hspace{-1cm}=
  \frac{1}{\sqrt{\pi}}\,\frac{1}{2^{\ms}\,\ms!}
  \,\frac{\asp}{\alpha}\sum_{\ell=0}^{N}\ush_{\ell}(\ts)
  \bigg[
    \int_{\REAL}\Hs_{\ell}(\zeta)\Hs_{\ms}(\zeta)e^{-\zeta^2}\,d\zeta +
    \int_{\REAL}\Hs_{\ell}(\zeta)\,\zeta\Hsp_{\ms}(\zeta)e^{-\zeta^2}\,d\zeta
    \bigg]
  \nonumber\\[0.5em]
  &\qquad
  = \TERM{I}{} + \TERM{II}{}.
\end{align}
These are finally evaluated as follows:
\begin{align}
  \TERM{I}{}
  &=
  \frac{1}{\sqrt{\pi}}\,\frac{1}{2^{\ms}\,\ms!}
  \frac{\asp}{\alpha}
  \sum_{\ell=0}^{N}\ush_{\ell}(\ts)\big(2^{\ms}\,\ms!\sqrt{\pi}\big)\delta_{\ell,\ms}
  =
  \frac{\asp(\ts)}{\alpha}\ush_{\ms}.
  \nonumber\\[1.25em]
  %% --------------------------------------------------------------------------------------------------------------------
  \TERM{II}{}
  &=
  \frac{1}{\sqrt{\pi}}\,\frac{1}{2^{\ms}\,\ms!}
  \,\frac{\asp}{\alpha}\sum_{\ell=0}^{N}\ush_{\ell}(\ts)
  \bigg[
    \int_{\REAL}\Hs_{\ell}(\zeta)\,\zeta\Hsp_{\ms}(\zeta)e^{-\zeta^2}\,d\zeta
    \bigg]
  \nonumber\\[0.5em]
  &=  
  \frac{1}{\sqrt{\pi}}\,\frac{1}{2^{\ms}\,\ms!}
  \frac{\asp}{\alpha}
  \sum_{\ell=0}^{N}\ush_{\ell}(\ts)
  \big(2^{\ell}\,\ell!\sqrt{\pi}\big)
  \bigg[ 2\ms(\ms-1)\delta_{\ms-2,\ell} + \ms\delta_{\ell,\ms} \bigg]
  \nonumber\\[0.5em]
  &=
  \frac{1}{\sqrt{\pi}}\,\frac{1}{2^{\ms}\,\ms!}
  \frac{\asp}{\alpha}
  \big(2^{\ms-2}\,(\ms-2)!\sqrt{\pi}\big)
  \bigg[ 2\ms(\ms-1) \ush_{\ms-2}(\ts) + \big(2^{\ms}\,\ms!\sqrt{\pi}\big)\,\ms\ush_{\ms}(\ts) \bigg]
  \nonumber\\[0.5em]
  &=
  \frac{\asp}{\alpha}\bigg[
    \frac12\ush_{\ms-2}(\ts) + \ms\ush_{\ms}(\ts)
    \bigg].
\end{align}
The final result of all these computations is again the scheme
\eqref{eq:form0:final:15}.

\end{document}